\documentclass[12pt]{amsart}

\usepackage{fullpage}
\usepackage{hyperref}
\usepackage{amssymb,amsthm,amsmath,amsrefs,amsfonts}
\usepackage{{latexsym}}

\numberwithin{equation}{section}

\newtheorem{theorem}{Theorem}[section]
\newtheorem{lemma}[theorem]{Lemma}
\newtheorem{corollary}[theorem]{Corollary}
\newtheorem{proposition}[theorem]{Proposition}

\theoremstyle{remark}
\newtheorem{remark}[theorem]{Remark}

\theoremstyle{definition}

\newcommand{\QT}{\mathrm{QT}}
\newcommand{\F}{\mathrm{F}}
\newcommand{\LL}{\mathrm{L}}
\newcommand{\Cu}{\mathrm{Cu}}
\newcommand{\T}{\mathrm{T}}
\newcommand{\ped}{\mathrm{Ped}}
\newcommand{\lat}{\mathrm{Lat}}

\newcommand{\LS}{\mathrm{Lsc}}

\newcommand{\mcCu}{\mathcal{C}u}
\newcommand{\mcC}{\mathcal{C}}

\newcommand{\N}{\mathbb{N}}

\newcommand{\R}{\mathbb{R}}

\newcommand{\wayb}{\ll}

\newcommand{\Ideal}{\mathrm{Ideal}}
\newcommand{\Set}{\mathrm{Set}}

\DeclareMathOperator{\rank}{rank}
\DeclareMathOperator{\supp}{fin}

\title[The cone of lower semicontinuous traces on a C*-algebra]
{The cone of lower semicontinuous traces on a C*-algebra}

\author{George. A. Elliott \and Leonel Robert \and Luis Santiago}

\address{George A. Elliott, Department of Mathematics, University of
 Toronto, Toronto, Canada~ M5S 2E4}
\email{elliott@math.toronto.edu}
\address{Leonel Robert, The Fields Institute, Toronto, Canada~ M5T 3J1}
\email{lrobert@math.toronto.edu}
\address{Luis Santiago, Department of Mathematics, University of
 Toronto, Toronto, Canada~ M5S 2E4}
\email{santiago@math.toronto.edu}

\begin{document}

\begin{abstract}
The cone of lower semicontinuous traces is studied with a view to
its use as an invariant. 
Its properties include compactness, Hausdorffness, and
continuity with respect to inductive limits. A suitable notion of dual cone is given. 
The cone of lower semicontinuous 2-quasitraces on a
(non-exact) C*-algebra is considered as well.
These results are applied to the study of the Cuntz semigroup. 
It is shown that if a C*-algebra absorbs the Jiang-Su algebra, then 
the subsemigroup of its Cuntz semigroup consisting of the purely non-compact elements is isomorphic to the dual
cone of the cone of lower semicontinuous 2-quasitraces. This yields a computation
of the Cuntz semigroup for the following two classes of C*-algebras: C*-algebras that absorb 
the Jiang-Su algebra and have no non-zero simple subquotients, and simple C*-algebras that absorb the 
Jiang-Su algebra.
\end{abstract}

\maketitle

\section{Introduction}

The most standard invariants in the classification of nuclear, simple, C*-algebras are their K-groups and their traces.
The traces are assumed to be bounded in the unital case, and lower semicontinuous and densely finite in the non-unital case.
If one has in mind the classification of non-simple C*-algebras, it is clear that these two kinds of 
traces will not suffice and a 
broader class should be considered. In this paper we study the properties of the cone of all lower semicontinuous
traces on a C*-algebra, with the purpose of applying our results to questions in the classification of 
non-simple C*-algebras.
 We also consider lower semicontinuous 2-quasitraces, since they appear naturally as functionals on the Cuntz
semigroup of the algebra. If the algebra is exact (this is the case that we are mostly concerned with in
the classification program), then lower semicontinuous traces and 2-quasitraces coincide. However, some 
of our considerations apply equally to traces and 2-quasitraces without assuming exactness of the C*-algebra.
Thus, we treat both classes for arbitrary C*-algebras.

Recall that a trace on a C*-algebra $A$ is a linear map $\tau$ on the positive
elements of $A$, with values in $[0,\infty]$,
that vanishes at 0 and satisfies the trace identity 
$\tau(xx^*)=\tau(x^*x)$ (see \cite{dixmier}). Every trace on $A$ extends
to a trace on $A\otimes \mathcal K$. A 2-quasitrace on $A$ is a map
on $(A\otimes \mathcal K)^+$, 
with values in $[0,\infty]$, that vanishes at 0,
satisfies the trace identity, and is linear on pairs of positive elements
that commute (see \cite[Definition 2.22]{kirchberg-blanchard} and
\cite[Proposition 2.24]{kirchberg-blanchard}). 
If a trace or 2-quasitrace is lower semicontinuous, then it is invariant 
under approximately inner automorphisms. 
This makes the cones of lower semicontinuous traces and 2-quasitraces
on a C*-algebra $A$---let us denote these cones by 
$\T(A)$ and $\QT_2(A)$---natural classification invariants associated to $A$.

In Section 2 of this paper we study the basic properties of the cone $\T(A)$.
Even though various classes of traces have been studied in the past 
(e.g., \cite{dixmier}, \cite{pedersen}, \cite{perdrizet}), 
we have found no bibliographic source for the properties of $\T(A)$. On the other hand,
the subcone of $\T(A)$ consisting of densely finite traces has been studied more thoroughly 
(e.g., in \cite{pedersen} and \cite{perdrizet}). 
We will show that some well known properties of the cone of densely finite traces
persist as properties of $\T(A)$. Furthermore, some properties appear that are not present in the cone
of densely finite traces, notably, the compactness of $\T(A)$ (in a suitable topology). 

In Section 3 we turn our attention to $\QT_2(A)$. By results of Blanchard and Kirchberg (see 
\cite[Proposition 2.24]{kirchberg-blanchard}),
which in turn extend work by Cuntz, Blackadar, Handelman, and Goodearl, 
the lower semicontinuous 2-quasitraces are in bijective correspondence
with the additive, order-preserving, extended positive real-valued maps on the Cuntz semigroup that vanish at 0 and preserve the suprema of 
increasing sequences---which henceforth we shall just call
functionals. Thus, we may think of $\QT_2(A)$ as the cone of functionals on the Cuntz semigroup.

Section 4 contains the description of suitable dual cones for $\QT_2(A)$ and $\T(A)$. The main results
of this section, Theorem \ref{supsofas} and Theorem \ref{supsofasT}, relate the dual cones
of $\QT_2(A)$ and $\T(A)$ with the functions that the positive elements of $A$ induce 
on $\QT_2(A)$ and $\T(A)$.

The last section contains applications of our results to understanding the structure of the Cuntz semigroup of 
certain C*-algebras,
in particular those C*-algebras that absorb the Jiang-Su algebra. For this class, we identify a natural subsemigroup
of the Cuntz semigroup that is isomorphic to the dual cone of $\QT_2(A)$. The complement of this subsemigroup
consists of the elements that become  compact, and not a multiple of infinity, after passing to the quotient 
by some closed two-sided ideal. The last result
of the paper is the computation of the Cuntz semigroup for two (disjoint) classes of C*-algebras: 
C*-algebras that absorb the Jiang-Su algebra and have no 
non-zero simple subquotients, and simple C*-algebras that absorb the Jiang-Su algebra.
The computation of the Cuntz semigroup for the latter class extends a previous result of Brown, Perera, and Toms 
(see \cite{brown-perera-toms}); 
in their computation they made the additional assumptions that the algebra was unital, exact, and of stable rank one.

\section{Preliminary results}
\subsection{The Cuntz-Pedersen equivalence relation}
Let $A$ be a C*-algebra and  let $a$ and $b$ be positive elements of $A$. Let us say that $a$
is Cuntz-Pedersen equivalent to $b$, and write $a\sim b$, if $a=\sum_{i=1}^\infty x_ix_i^*$ and 
$b=\sum_{i=1}^\infty x_i^*x_i$
for some sequence $x_i\in A$, $i=1,2,\dots$. Let us say that $a$ is Cuntz-Pedersen smaller than $b$
if $a\sim a'$ for some $a'\in A^+$ with $a'\leq b$. In this case let us write $a\preccurlyeq b$. In 
\cite{cuntz-pedersen} Cuntz and Pedersen showed 
that the relations $\sim$ and $\preccurlyeq$ are transitive.  
This is a consequence of the following result of Pedersen 
(\cite[Proposition 1.1]{pedersen}).

\begin{proposition} (Riesz-Pedersen decomposition property.)
Suppose that $x_i,y_i\in A$, $i=1,2\dots$,  are such that 
$\sum_{i=1}^\infty x_ix_i^*=\sum_{i=1}^\infty y_i^*y_i$. 
Then there are elements $z_{i,j}$, $i,j=1,2\dots$, such that $x_i^*x_i=\sum_{j=1}^\infty z_{i,j}^*z_{i,j}$ and 
$y_jy_j^*=\sum_{i=1}^\infty z_{i,j}z_{i,j}^*$, for all $i,j\geq 1$.
\end{proposition}

\begin{lemma}\label{K-R}
Let $a,b\in A^+$ be such that $\|a-b\|<\epsilon$. Then 
$(a-\epsilon)_+\preccurlyeq b$.
\end{lemma}
\begin{proof}
By \cite[Lemma 2.2]{kirchberg-rordam}, there is $d\in A$ with $\|d\|\leq 1$ and
$(a-\epsilon)_+=dbd^*$. Hence, $(a-\epsilon)_+\sim b^{1/2}d^*db^{1/2}\leq b$.
\end{proof}

The following proposition is a summary of the properties of the relations $\preccurlyeq$
and $\sim$ between positive elements of $A$ that will be needed later.

\begin{proposition}\label{CP}
(i) For every $\epsilon>0$ there is $\epsilon'>0$
such that
\begin{eqnarray}
(a_1+a_2-\epsilon)_+ &\preccurlyeq& (a_1-\epsilon')_++(a_2-\epsilon')_+ \label{CP1},\\
(a_1-\epsilon)_++(a_2-\epsilon)_+ &\preccurlyeq & (a_1+a_2-\epsilon')_+ \label{CP2}.
\end{eqnarray}

(ii) For every $\epsilon>0$ and $x\in A$ we have $(xx^*-\epsilon)_+\sim (x^*x-\epsilon)_+$.
 
(iii) (Riesz decomposition property.)  Suppose that $\sum_{i=1}^\infty a_i\sim \sum_{i=1}^\infty b_i$. 
Then there are elements $w_{i,j}$, $i,j=1,2,\dots$, such that 
$a_i\sim \sum_{j=1}^\infty w_{i,j}$ and $b_j\sim \sum_{i=1}^\infty w_{i,j}$.
\end{proposition}

\begin{proof}
(i) Since $(a_1-\epsilon')_+ + (a_2-\epsilon')_+$ tends to $a_1+a_2$
as $\epsilon'$ tends to 0, we get \eqref{CP1} from Lemma \ref{K-R}. 

Let us prove \eqref{CP2}. Let $(e_n)$ be an approximate unit for $\overline{(a_1+a_2)A(a_1+a_2)}$ such that 
$e_n(a_1+a_2)e_n\leq (a_1+a_2-1/n)_+$ (e.g., $e_n=\phi_n(a_1+a_2)$, with $\phi_n(t)=\frac 1 t(t-1/n)_+$).
Since $e_na_1e_n\to a_1$ and $e_na_2e_n\to a_2$, by Lemma \ref{K-R} there exists $n$ such that
$(a_1-\epsilon)_+ + (a_2-\epsilon)_+\preccurlyeq e_n(a_1+a_2)e_n\leq (a_1+a_2-1/n)_+$.

(ii) Let us show that $(x^*x-\epsilon)_+$ is Murray-von Neumann equivalent to $(xx^*-\epsilon)_+$
(i.e., Cuntz-Pedersen equivalent by means of a single element). 
Consider the polar decomposition $x=u|x|$ of $x$ in the bidual of $A$. 
The element $y=u(x^*x-\epsilon)_+^{1/2}$ belongs to $A$, and we have $y^*y=(x^*x-\epsilon)_+$ and 
$yy^*=(xx^*-\epsilon)_+$.

(iii) Suppose that $\sum_i a_i=\sum_j x_j^*x_j$ and $\sum_k b_k=\sum_j x_jx_j^*$. By the Riesz-Pedersen
decomposition property, there are $y_{i,j}$s such that $a_i=\sum_j y_{i,j}^*y_{i,j}$ and 
$x_jx_j^*=\sum_i y_{i,j}y_{i,j}^*$. We have $\sum_{i,j} y_{i,j}y_{i,j}^*=\sum_k b_k$. Therefore,
there are $z_{i,j,k}$s such that $b_k=\sum_{i,j} z_{i,j,k}^*z_{i,j,k}$ and 
$y_{i,j}^*y_{i,j}=\sum_k z_{i,j,k}z_{i,j,k}^*$. 
Set $\sum_j z_{i,j,k}^*z_{i,j,k}=w_{i,k}$. Then 
$\sum_k w_{i,k}=\sum_{j,k} z_{i,j,k}^*z_{i,j,k}\sim \sum_j y_{i,j}^*y_{i,j}=a_i$. Also, 
$\sum_i w_{i,k}\sim\sum_{i,j} z_{i,j,k}z_{i,j,k}^*=b_k$. 
\end{proof}

\subsection{Non-cancellative cones.}
Let us introduce the terminology  non-cancellative cone for an abelian  semigroup endowed with a scalar 
multiplication by strictly positive real numbers. 
The semigroup may not have cancellation, that is to say, 
$\tau+\tau_1=\tau+\tau_2$ may not imply that $\tau_1=\tau_2$. However, we will often 
refer to non-cancellative cones simply as cones, and
refer to standard cones that embed in a vector space as cancellative cones.

Notice that we have not included scalar multiplication by 0 or $\infty$ in the definition
of non-cancellative cone. For some of the cones that we shall consider here---of traces and, more generally, of 2-quasitraces---we 
will be able to extend the scalar multiplication to include $0$ and $\infty$. However, it will
not  necessarily be the case that scalar multiplication by $0$ will result in the zero element 
of the cone (see \eqref{extendedsc}).

Non-cancellative cones satisfy the following form of restricted cancellation. 
 
\begin{lemma} \label{cancellation}
(Cancellation lemma.) Let $S$ be a non-cancellative cone. 
Suppose that $x+z=y+z$ for some $z$ such that $z+z_1=nx$ and $z+z_2=ny$. Then
$x=y$.
\end{lemma}
\begin{proof}
By induction we have $nx+z=ny+z$. So let us assume that $n=1$. Then $x+y=x+z+z_2=y+z+z_2=2y$.
In the same way $x+y=2x$, and so $x=y$.
\end{proof}

\section{The cone of lower semicontinuous traces}

\subsection{The cone $\T(A)$}
Let $A$ be a C*-algebra. Let us say that $\tau\colon A^+\to [0,\infty]$
is a trace on $A$ if $\tau$ is linear (i.e., additive, homogeneous 
with respect to strictly
positive scalars, and vanishing at 0) and satisfies the trace 
identity $\tau(xx^*)=\tau(x^*x)$.
 
The following lemma is well known (see \cite[Remark 2.27 (iv)]{kirchberg-blanchard}).

\begin{lemma}\label{whatsatrace}
If $\tau\colon A^+\to [0,\infty]$ is a trace  then $\widetilde\tau$ defined by
$\widetilde\tau(a):=\sup_{\epsilon>0}\tau((a-\epsilon)_+)$ is a 
lower semicontinuous trace, and is the largest such trace
majorized by $\tau$.

\end{lemma}
\begin{proof}
By parts (i) and (ii) of Proposition \ref{CP} we have that
$\widetilde\tau$ is a trace.
We also have that $\widetilde\tau(a)=\sup_{\epsilon>0}\widetilde\tau((a-\epsilon)_+)$. 
Let us show that this implies that $\widetilde\tau$ is lower
semicontinuous. Suppose that $\widetilde\tau(a)>\alpha$ 
for some $\alpha\geq 0$.
Let $\epsilon>0$ be such that $\widetilde\tau((a-\epsilon)_+)>\alpha$. By Lemma \ref{K-R}, if $\|a'-a\|<\epsilon$ 
then $(a-\epsilon)_+\preccurlyeq a'$, whence 
$\widetilde\tau(a')\geq \widetilde\tau((a-\epsilon)_+)>\alpha$.
If $\sigma$ is another lower semicontinuous trace
with $\sigma\le \tau$, then for any $a\in A^+$,

\[
\sigma (a)=
\sup_{\epsilon>0} \sigma((a-\epsilon)_+)\le \sup_{\epsilon>0}
\tau ((a-\epsilon)_+)=\widetilde\tau (a)
\qedhere
\]
\end{proof}

Let us denote by $\T(A)$ the collection of all lower semicontinuous traces of $A$. This set is
a non-cancellative cone endowed with the operations of pointwise addition and pointwise scalar multiplication by 
strictly positive real numbers. (We will later extend the scalar multiplication to include $0$ and $\infty$.) 
We shall also consider $\T(A)$ endowed with the order induced by its addition operation. (When we consider
the dual cone of $\T(A)$ in Section \ref{duals} below, we shall also need to consider its pointwise order, but, as we shall now show, for $\T(A)$ itself this is determined by addition.) 

The following proposition is well known for various 
classes of traces on a C*-algebra (e.g., see \cite[Proposition 6]{dixmier}).
\begin{proposition}\label{orderoftraces}
Let $\tau_1,\tau_2\in \T(A)$. Suppose that $\tau_1(x)\leq \tau_2(x)$ for all
$x\in A^+$. Then  there exists $\tau_3\in \T(A)$ such that $\tau_1+\tau_3=\tau_2$. 
\end{proposition}
\begin{proof}
Define $\tau\colon A^+\to [0,\infty]$ as follows:
\[
\tau(x) :=\left \{
\begin{array}{cl}
\tau_2(x)-\tau_1(x) &\hbox{  if }\tau_2(x)<\infty, \\
\infty & \hbox{ otherwise.}\\
\end{array}\right.
\]
It is easily verified that $\tau$ is linear, satisfies the trace identity, and  satisfies $\tau_1+\tau=\tau_2$.
Set $\widetilde\tau=\tau_3$, where $\widetilde\tau$ is the lower semicontinuous regularization of $\tau$ described in 
Lemma \ref{whatsatrace}. Taking the suprema of both sides with respect to $\epsilon$ in the 
equation 
\[
\tau_1((a-\epsilon)_+)+\tau((a-\epsilon)_+)=\tau_2((a-\epsilon)_+),
\]
we get that $\tau_1+\tau_3=\tau_2$.
\end{proof}

In \cite[Theorem 3.1]{pedersen}, Pedersen used the Riesz-Pedersen property to show that the cone of densely finite 
lower semicontinuous traces is a lattice. We shall follow a similar method here to show that the whole of
$\T(A)$ is a lattice, and is in fact complete.

\begin{theorem}\label{lattice}
The cone $\T(A)$ is a complete lattice with respect to the order determined by addition
(equivalently, by Proposition \ref{orderoftraces}, the pointwise order). For all $\tau_1,\tau_2,\tau_3\in \T(A)$
we have  
\begin{eqnarray}\label{supidentity}
\tau_1\vee\tau_2+\tau_3=(\tau_1+\tau_3)\vee (\tau_2+\tau_3),\\
\tau_1\wedge\tau_2+\tau_3=(\tau_1+\tau_3)\wedge (\tau_2+\tau_3).\label{infidentity} 
\end{eqnarray}
\end{theorem}

\begin{proof}
The properties of lower semicontinuity and  linearity, and also the trace identity, are preserved under passing
to the pointwise supremum of an upward directed collection of lower semicontinuous traces. 
Thus, $\T(A)$ is closed under passage to  directed suprema. In order to prove that
$\T(A)$ is a complete lattice it is then enough to show that
the supremum of any two lower semicontinuous traces exists. (The supremum of any non-empty
set will then exist, and the supremum of the empty set is 0. It follows that the infimum of 
any set also exists.) 

Let $\tau_1$ and $\tau_2$ be in $\T(A)$. Define $\tau\colon A^+\to[0,\infty]$ by  the Riesz-Kantorovich formula: 
\[
\tau(x):=\sup \{\,\tau_1(x_1)+\tau_2(x_2)\mid x_1+x_2\sim x\,\}.
\]
We clearly have $\tau(xx^*)=\tau(x^*x)$. The linearity of $\tau$ follows
from the Riesz decomposition property (i.e., Proposition \ref{CP} (iii)), by a standard argument that goes back to Riesz
(see \cite[Theorem 1]{riesz}). It is clear that $\tau_1\leq \tau$, $\tau_2\leq \tau$, and that any trace
that majorizes $\tau_1$ and $\tau_2$ is greater than or equal to $\tau$. 
With $\widetilde\tau$ the lower semicontinuous regularization of $\tau$ of Lemma \ref{whatsatrace}, i.e., 
$\widetilde\tau(a)=\sup_{\epsilon>0} \tau((a-\epsilon)_+)$, we have  
$\widetilde \tau\leq \tau$, and, furthermore, both  $\tau_1\leq \widetilde\tau$ and $\tau_2\leq \widetilde\tau$
(as $\tau_i=\sup_{\epsilon>0}\tau_i((a-\epsilon)_+)\leq\sup_{\epsilon>0}
\tau((a-\epsilon)_+)=\widetilde\tau(a)$). Therefore,
$\tau\leq \widetilde\tau$, and so $\tau=\widetilde\tau$; in other 
words, the supremum of $\tau_1$ and $\tau_2$ in the cone of all traces belongs to $\T(A)$.

The identity \eqref{supidentity} follows from the Riesz-Kantorovich formula for the supremum
of two traces in $\T(A)$. (Note that we have shown that this formula does describe the supremum
in $\T(A)$.)

Consider now $\tau\colon A^+\to [0,\infty]$ defined by
\[
\tau(x):=\inf\{\,  \tau_1(x_1)+\tau_2(x_2)\mid x_1+x_2\sim x\,\}.
\]
By the Riesz decomposition property $\tau$ is a trace. Very much as shown above for the
supremum, $\tau$ is seen to be lower semicontinuous and to be the infimum of $\tau_1$ and
$\tau_2$ in $\T(A)$. The identity  \eqref{infidentity} now follows from
the Riesz-Kantorovich formula that we have just established for the infimum of $\tau_1$ and
$\tau_2$ in $\T(A)$.
\end{proof}

Vector lattices, i.e., ordered vector spaces that are lattices with respect to their order,
have a number of properties that are implied by their lattice structure. For example,  
a vector lattice  is always distributive and satisfies the
identities \eqref{supidentity} and \eqref{infidentity} (see \cite{vector-lattices}). The cone $\T(A)$ 
cannot be embedded in a vector space
since it is not cancellative. For instance, if $I$ denotes a closed two-sided
ideal of $A$ then $\tau_I$ defined by
\[
\tau_I(x):=\left\{
\begin{array}{lc}
0 & x\in I^+,\\
\infty & x\notin I^+,
\end{array}\right. 
\]
is a lower semicontinuous trace and satisfies $\tau_I+\tau_I=\tau_I$. 
Indeed, the lower semicontinuous traces 
with the only possible values $0$ and 
$\infty$---i.e., that satisfy $\tau+\tau=\tau$---are, as is easily seen, in order reversing bijection with the 
closed two-sided ideals of $A$ by the map $I\mapsto \tau_I$. 

Making use of equations \eqref{supidentity} and \eqref{infidentity},
and the restricted cancellation of Lemma \ref{cancellation}, we can show that $\T(A)$ has 
some of the properties of a vector lattice.

\begin{proposition}\label{distributive}
(i) We have $\tau_1\vee\tau_2+\tau_1\wedge\tau_2=\tau_1+\tau_2$ for all $\tau_1$ and $\tau_2$  in $\T(A)$. 

(ii) $\T(A)$ is a distributive lattice.
\end{proposition} 

\begin{proof}
(i) Taking $\tau_3=\tau_1\wedge \tau_2$ in \eqref{supidentity} yields 
\[\tau_1\vee\tau_2+\tau_1\wedge\tau_2=
(\tau_1+\tau_1\wedge \tau_2)\vee (\tau_2+\tau_1\wedge \tau_2)\leq \tau_1+\tau_2.
\] 
Taking $\tau_3=\tau_1\vee \tau_2$ in \eqref{infidentity} yields  
\[
\tau_1\wedge\tau_2+\tau_1\vee\tau_2=
(\tau_1+\tau_1\vee \tau_2)\wedge (\tau_2+\tau_1\vee \tau_2)\geq \tau_1+\tau_2.
\]

(ii) Let us prove that $(\tau_1\vee\tau_2)\wedge\tau_3=(\tau_1\wedge\tau_3)\vee(\tau_2\wedge\tau_3)$.
It is enough to prove this equality after adding  
$\tau_1\wedge\tau_2\wedge\tau_3$ to both sides, since this term may be cancelled by Lemma \ref{cancellation}. 
Considering the right-hand side, we have
\begin{align*}
(\tau_1\wedge\tau_3)\vee(\tau_2\wedge\tau_3)+\tau_1\wedge\tau_2\wedge\tau_3 &=
\tau_1\wedge\tau_3+\tau_2\wedge\tau_3\\
&=(\tau_1+\tau_2\wedge\tau_3)\wedge(\tau_3+\tau_2\wedge\tau_3) \\
&=(\tau_1+\tau_2)\wedge(\tau_1+\tau_3)\wedge(\tau_2+\tau_3)\wedge 2\tau_3.
\end{align*}
Considering the left-hand side, we obtain the same quantity:
\begin{align*}
(\tau_1\vee\tau_2)\wedge\tau_3+\tau_1\wedge\tau_2\wedge\tau_3 
&=(\tau_1\vee\tau_2+\tau_1\wedge\tau_2\wedge\tau_3)\wedge (\tau_3+\tau_1\wedge\tau_2\wedge\tau_3)\\
&=(\tau_1+\tau_2)\wedge (\tau_1\vee\tau_2+\tau_3)\wedge (\tau_3+\tau_1\wedge\tau_2\wedge\tau_3)\\
&=(\tau_1+\tau_2)\wedge (\tau_3+\tau_1\wedge\tau_2\wedge\tau_3) \\
&=(\tau_1+\tau_2)\wedge(\tau_1+\tau_3)\wedge(\tau_2+\tau_3)\wedge 2\tau_3.
\qedhere
\end{align*}
\end{proof}

\subsection{The topology on $\T(A)$}
Let us endow the cone $\T(A)$  with the topology in which the net $(\tau_i)$ converges to $\tau$  if
\begin{align}\label{topology}
\limsup \tau_i((a-\epsilon)_+)\leq \tau(a)\leq \liminf \tau_i(a)
\end{align}
for any $a\in A^+$ and $\epsilon>0$. Equivalently 
(by Lemma \ref{whatsatrace}---both parts---combined with compactness of the 
infinite product of copies of $[0,\infty]$---one for each $a\in A^+$), $\tau_i\to \tau$ if,
whenever a subnet of $(\tau_i)$ converges pointwise to a trace $\sigma$, the regularization $\widetilde\sigma$
of $\sigma$ given by Lemma \ref{whatsatrace} is equal to $\tau$ (cf.~proof of Theorem \ref{compact-hausdorff}, below).
A sub-basis of neighbourhoods for the  trace $\tau$ is given by the sets
\begin{align*}
U(\tau;a,\epsilon)& :=\{\,\tau'\in \T(A)\mid \tau((a-\epsilon)_+)\leq \tau'(a)+\epsilon \hbox{ or }\tau'(a)>\frac 1 \epsilon\,\},\\
V(\tau;a,\epsilon)& :=\{\,\tau'\in \T(A)\mid \tau'((a-\epsilon)_+)\leq \tau(a)+\epsilon\,\}.
\end{align*}

\begin{remark} \label{afterdef}
In order to define the topology of $\T(A)$ the element $a$  
can be restricted to vary in a dense subset $S$ of $A^+$ such that
$a\in S$ implies that $(a-1/n)_+\in S$ for all $n\geq 1$. 
Let us verify this. Let $S$ be such a set. 
Let $a\in A^+$ and $\epsilon>0$. Choose $a'\in S$ and $n\in \N$ such that
$\|a-(a'-2/n)_+\|<\epsilon$ and $\|a-a'\|<1/n$.
By Lemma \ref{K-R} (applied twice), we have
\[
(a-\epsilon)_+\preccurlyeq (a'-\frac 2 n)_+\preccurlyeq (a'-\frac 1 n)_+\preccurlyeq a.
\] 
Hence, $U(\tau;(a'-\frac 1 n)_+,\frac 1 n)\subseteq U(\tau;a,\epsilon)$ 
and 
$V(\tau;(a'-\frac 1 n)_+,\frac 1 n)\subseteq V(\tau;a,\epsilon)$.
\end{remark}

One can verify using \eqref{topology}
that $\alpha\tau\to \tau_{\ker \tau}$ when $\alpha\to \infty$. One also verifies
that $\alpha\tau\to \tau_{\supp \tau}$ when $\alpha\to 0$, where $\supp \tau$
is the closed two-sided ideal spanned by $\{\,x\in A^+\mid \tau(x)<\infty\,\}$.
(In the terminology of \cite{dixmier}, $\supp \tau$ is the closure of the ideal of definition
of $\tau$; we shall refer to this ideal as the (closed) ideal of finiteness of $\tau$.)
In view of these computations, we may extend by continuity the scalar multiplication 
in order to include the scalars 0 and $\infty$:
\begin{equation}\label{extendedsc}
0\cdot\tau=\tau_{\supp \tau}, \quad \infty\cdot\tau=\tau_{\ker \tau}.
\end{equation}

\begin{proposition}\label{addition} 
Addition and (extended) scalar multiplication in $\T(A)$ are jointly continuous. 
\end{proposition}
\begin{proof}
We will prove here that if $\tau_i\to \tau$ then $0\cdot \tau_i\to 0\cdot \tau$.
The other parts of the proposition are easily verified from the definition of the topology of $\T(A)$
by the inequalities \eqref{topology}. 

Let us show that $\limsup (0\cdot \tau_i)((a-\epsilon)_+)\leq (0\cdot\tau)(a)$ for all $a\in A^+$ and $\epsilon>0$.
If $(0\cdot \tau)(a)=\infty$ this is obvious. Suppose that
$(0\cdot \tau)(a)=0$. Then $a\in \supp \tau$, and so $\tau((a-\epsilon)_+)<\infty$ for all $\epsilon>0$.
Since $\tau_i\to \tau$, we have $\limsup \tau_i((a-\epsilon)_+)\leq \tau((a-\epsilon/2)_+)<\infty$.
Hence $(0\cdot \tau_i)((a-\epsilon)_+)=0$ for all $i\geq i_0$ for some $i_0$, and from this the desired
inequality follows. 

Let us show now that $(0\cdot \tau)(a)\leq \liminf (0\cdot \tau_i)(a)$ for all $a\in A^+$. 
If $\tau((a-\epsilon)_+)<\infty$ for all $\epsilon>0$ then $(0\cdot \tau)(a)=0$, and so
the desired inequality clearly holds.  
Suppose that $\tau((a-\epsilon)_+)=\infty$ for some $\epsilon>0$. Then $\tau_i((a-\epsilon)_+)=\infty$
for all $i\geq i_0$, for some $i_0$. Hence  $(0\cdot\tau_i)(a)=\infty$ for all $i\geq i_0$,
and from this the desired
inequality follows.
\end{proof}

\begin{theorem}\label{compact-hausdorff} $\T(A)$ is a compact Hausdorff space. 
If $A$ is separable then $\T(A)$ has a countable basis.
\end{theorem}
\begin{proof}
Let us show that $\T(A)$ is Hausdorff.
Let $\tau_1$ and $\tau_2$ be distinct points in $\T(A)$. Since either
$\tau_1\nleq \tau_2$ or $\tau_2\nleq \tau_1$,  we may suppose that we are in the first case. Then there
are $a\in A^+$ and $\epsilon>0$ such that $\tau_1((a-\epsilon)_+)\nleq \tau_2(a)+\epsilon$.
Let us choose $\epsilon>0$ such that $\tau_2(a)<2/\epsilon-\epsilon/2$ (this is possible since, necessarily,
$\tau_2(a)<\infty$).  Then the sets $U(\tau_1;(a-\epsilon/2)_+,\epsilon/2)$ and $V(\tau_2;a,\epsilon/2)$
are disjoint neighbourhoods of $\tau_1$ and $\tau_2$ respectively. 
For suppose that $\tau$ belongs to their intersection.
Then, either 
\[
\tau_1((a-\epsilon)_+)\leq \tau((a-\epsilon/2)_+)+\epsilon/2\leq \tau_2(a)+\epsilon
\] 
or 
\[
\frac 2 \epsilon<\tau((a-\epsilon/2)_+)\leq \tau_2(a)+\frac \epsilon 2.
\]
In either case, this is a contradiction. 
 
The following simple proof of the compactness of $\T(A)$ was suggested to 
us by E.~Kirchberg
(our original proof was much longer).

Let $(\tau_i)_{i\in \Lambda}$ be a net of traces in $\T(A)$. By Tychonoff's
theorem (using the compactness of $[0,\infty]$), 
we can choose a subnet $(\tau_i)_{i\in \Lambda'}$ converging pointwise
to $\sigma$.  The function $\sigma\colon A^+\to [0,\infty]$
is linear and satisfies the trace identity. With $\widetilde\sigma$ 
the lower semicontinuous trace 
of Lemma \ref{whatsatrace}, i.e.,
$\widetilde\sigma(a)=\sup_{\epsilon>0}\sigma((a-\epsilon)_+)$,
let us show that $(\tau_i)_{i\in \Lambda'}$ converges to
$\widetilde\sigma$ in $\T(A)$ (i.e., the inequalities \eqref{topology} 
are satisfied for all $a\in A^+$ and $\epsilon>0$). 
Let $a\in A^+$ and let $\epsilon>0$. Then
\[
\limsup_{i\in \Lambda'} \tau_i((a-\epsilon)_+)=\sigma((a-\epsilon)_+)\leq \widetilde\sigma(a)\leq \sigma(a)=\liminf_{i\in \Lambda'} \tau_i(a).
\]
(This convergence is also immediate from the alternative form of the definition.)
This shows that $\T(A)$ is compact. 

Now suppose that $A$ is separable. It follows from the remark made after the
definition of the topology of
$\T(A)$ that if $A$
is separable then $\T(A)$ is first countable,
and in fact there is a countable basis of symmetric
entourages for a uniform structure giving rise to the
topology of $\T(A)$; let us choose such a basis. Inspection of the
entourages described shows that not only are they symmetric but also the
corresponding neighbourhoods in the topology are open; we shall
assume therefore that our countable basis consists of such entourages.

Separability of $A$ also implies that there is
a countable dense subset of $\T(A)$: as a set of 
maps from $A^+$ to $[0, \infty]$, we may naturally identify
$\T(A)$ with a subset of $\Pi_{a\in  A^+} [0, \infty]$, 
and (since the maps in $\T(A)$ are lower semicontinuous)
in fact, as we shall now show, with a subset of
$\Pi_{a\in S} [0, \infty]$ where $S$ is a suitable countable
dense subset of $A^+$. Since lower semicontinuous functions 
are not determined on just any dense subset, we must choose
$S$ to consist of a countable dense subset of $A^+$ (any such subset) 
together with, for each $a\in S$, the set of all elements
$(a-\epsilon)^+$ with $\epsilon=1/n$, $n=1, 2, 3, \ldots$.
(It follows from Lemma \ref{K-R} that $S$ separates elements
of $\T(A)$: if $\tau, \tau' \in \T(A)$ and $\tau$ and $\tau'$
agree on $S$, then for any $a\in A^+$, with $b\in S$ such that
$\|a-b\| < \epsilon$, we have $(b-\epsilon)_{+} \preceq a$, in the
sense of Cuntz and Pedersen, and so
\[
\tau(b)=\sup_{n\ge 1} \tau((b-1/n)_+)=\sup_{n\ge 1}
\tau'((b-1/n)_+)\le \tau' (b).
\]
In other words, $\tau \le \tau'$, and so by symmetry 
$\tau =\tau'$.)

Not only is the embedding of $\T(A)$ in $\Pi_{a\in S} [0, \infty]$ injective, but also by the
alternative definition of the topology on $\T(A)$ the
inverse of this map, from the image with the coordinate-wise
topology, is continuous. In other words, as we shall
now show, if $\tau_i(a)\to \tau(a)$ for all  $a\in S$, with 
$\tau_i$ and $\tau$ in $\T(A)$, then $\tau_i \to \tau$ in $\T(A)$.
It is enough to show that if $\tau'$ is a trace and
$\tau_i \to \tau'$ pointwise on $A^+$ then $(\tau')^{\sim}=\tau$.
By hypothesis, $\tau'$ agrees with $\tau$ on
$S$. By the choice of $S$, $(\tau')^{\sim}$ also coincides
with $\tau$ on $S$, and therefore by injectivity
$(\tau')^{\sim}=\tau$ in $\T(A)$.
Hence, one obtains a countable dense subset of $\T(A)$ as
the image  under the inverse map of a countable dense subset of its
domain---which exists as the countable
Cartesian product is a metrizable compact space.

It follows that $\T(A)$ has a countable basis for the 
topology under consideration, namely, the collection
of all neighbourhoods of a fixed dense sequence
$\tau_1, \tau_2, \ldots$ in $\T(A)$ corresponding to the
countable basis of symmetric entourages for the
uniform structure referred to above. (The proof of this
is just as if the entourages were determined by a metric, as the degrees of closeness
corresponding to a sequence of distances converging to zero. Let $\tau$ be a
point in $\T(A)$, and let $W$ be an arbitrary open neighbourhood of a 
symmetric entourage $U$ such that if $(\tau, s)\in U$
then $s\in W$. Choose an entourage $V$ such that if $(\tau, \tau')\in V$ and
$(\tau', \sigma)\in V$ then $(\tau, \sigma)\in U$.
We may choose $V$ to be one of the countable basis of symmetric entourages chosen 
above and in particular  such that the neighbourhood of
any point determined by $V$ is open.
Choose $n$ such that $(\tau, \tau_n)\in V$.
The neighbourhood of $\tau_n$ determined by the symmetric
entourage $V$ then both includes the point $\tau$ and
is included in the neighbourhood of $\tau$
determined by $U$, and therefore also in the
given open neighbourhood $W$ of $\tau$. Since this  
neighbourhood is open by the choice of $V$, we
have identified a countable basis of open sets
(the neighbourhoods of $\tau_1, \tau_2, \cdots$
determined by the chosen countable basis of entourages).)
\end{proof}

\begin{proposition}\label{infdistributivity}
(i) The order relation in $\T(A)$ is continuous (i.e., the set $\{(\tau_1,\tau_2)\mid \tau_1\leq \tau_2\}$ is closed
in $\T(A)\times \T(A)$).

(ii) An upward directed subset of $\T(A)$ converges to its supremum (when indexed by itself), and  a downward  directed subset converges 
to its infimum.

(iii) The complete, distributive, lattice $\T(A)$ is join continuous; that is, for any subset $S$ of  $\T(A)$, and for any
$\tau\in \T(A)$,
\[
(\bigwedge_{\sigma\in S}\sigma)\vee \tau=\bigwedge_{\sigma\in S}(\sigma\vee\tau).
\]
\end{proposition}
\begin{proof}
(i) Let $((\tau_i,\mu_i))$ be a net converging to $(\tau,\mu)$ and suppose that
$\tau_i\leq \mu_i$ for all $i$. Then $\mu_i=\tau_i+\mu_i'$ for some $\mu_i'$. Passing
to a convergent subnet of $\mu_i'$ (by compactness) and then passing to the limit we get that 
$\mu=\tau+\mu'$, whence $\tau\leq \mu$.

(ii) Let $(\tau_i)_{i\in \Lambda}$ be a decreasing net with infimum $\tau$. It is enough by compactness to show that
every convergent subnet of $(\tau_i)_{i\in \Lambda}$ converges to $\tau$, and so we may assume  without loss of generality that
$(\tau_i)_{i\in \Lambda}$ converges to $\tau'$. For every $i$ we have $\tau\leq \tau_i$. Thus, passing to the limit and using part (i)
of this proposition we conclude that $\tau\leq \tau'$. On the other hand, for every $i$ and $j$ with $i\leq j$
we have $\tau_j\leq \tau_i$. Fixing $i$ and passing to the limit in $j$ we obtain $\tau'\leq \tau_i$.
Since this holds for all $i$ we conclude that $\tau'\leq \tau$. 

One may proceed in a similar way for upward directed subsets of $\T(A)$.

(iii) By the distributivity of $\T(A)$ we have 
\[
(\bigwedge_{i\in F} \tau_i)\vee \tau=\bigwedge_{i\in F} (\tau_i\vee\tau)\] 
for every finite subset $F$ of $S$. Let us consider both sides as downward directed families of traces indexed by the finite
subsets of $S$. The infimum of the right side is $\bigwedge_{i\in \Lambda}(\tau_i\vee\tau)$. Set $\bigwedge_{i\in F} \tau_i=\mu_F$.
It is enough to prove that if  the  downward directed subset $\{\mu_F\}$ has infimum $\mu$, then the infimum 
of $\{\mu_F\vee\tau\}$ is $\mu\vee \tau$. By Proposition \ref{distributive} (i), we have $\mu_F\vee\tau+\mu_F\wedge\tau=\mu_F+\tau$.
By (ii) together with Proposition \ref{addition}, for any two downward directed sets $S_1$ and $S_2$,
$\bigwedge (S_1+S_2)=\bigwedge S_1+\bigwedge S_2$. Hence,
taking infima  on both sides of $\mu_F\vee\tau+\mu_F\wedge\tau=\mu_F+\tau$ we get 
\[\bigwedge_F (\mu_F\vee\tau)+\bigwedge_F(\mu_F\wedge \tau)=  
\bigwedge_F\mu_F+\tau=(\bigwedge_F\mu_F)\wedge\tau+(\bigwedge_F\mu_F)\vee\tau.\] 
Cancelling $(\bigwedge_F\mu_F)\wedge\tau$  (using Lemma \ref{cancellation} with $n=1$) we obtain
$\bigwedge_F (\mu_F\vee\tau)=(\bigwedge_F\mu_F)\vee\tau$, as desired.
\end{proof}

\begin{remark} Proposition \ref{infdistributivity} (i) and (ii) may be proved directly from the definition of the topology of
$\T(A)$. The proof given above, however, applies to an arbitrary topological cone that is a complete lattice, and 
is compact and Hausdorff. The infinite distributivity of Proposition \ref{infdistributivity} (iii) implies that the lattice 
obtained by reversing the order of $\T(A)$
is a continuous lattice (in the sense of \cite{compendium}; see \cite[Theorem I-2.7]{compendium}). Since the map $I\mapsto \tau_I$ is an order reversing embedding
of $\lat(A)$ as a subcomplete sublattice of $\T(A)$, we deduce from Proposition \ref{infdistributivity} (iii) the well
known fact that $\lat(A)$ is a continuous lattice.
\end{remark}

\emph{Question.} Is the map $\mu\mapsto \mu\vee \tau$ continuous in the topology of $\T(A)$? Is the topology of $\T(A)$ the Fell-Lawson topology of 
the complete lattice obtained by reversing the order of $\T(A)$?

Let us write $\T_I(A)$ for the subcone of $\T(A)$ of traces with ideal of finiteness $I$ (i.e., $\tau$ with $\supp \tau=I$).
The subcone $\T_I(A)$ is in bijective correspondence with the densely finite lower semicontinuous traces
on the ideal $I$, because for every such trace we get a trace in $\T_I(A)$ by setting it equal to $\infty$
outside $I^+$.

\begin{proposition}\label{inducedtops}
(i) For each ideal $I$ of $A$, the relative topology on the subcone $\T_I(A)$ of $\T(A)$
is the topology of pointwise convergence on the positive elements of the Pedersen ideal of $I$. 

(ii) The relative topology on the subset $\lat(A)$ of $\T(A)$---the image of the embedding
$I\mapsto \tau_I$---is the Fell-Lawson topology.
\end{proposition}
\begin{proof}
(i) Let $\tau\in \T_I(A)$ and let $(\tau_i)$ be a net in $\T_I(A)$ converging to $\tau$ in the topology of $\T(A)$. 
Let us show that $\tau_i(a)\to \tau(a)$ for all $a\in \ped(I)^+$. By the alternative definition of limit
in $\T(A)$ as the regularization of every pointwise convergent subnet, it is sufficient to show that an arbitrary
densely finite trace $\sigma$ on $I^+$ satisfies $\sigma=\sup_{\epsilon>0} \sigma((a-\epsilon)_+)$ for each
$a\in \ped(I)^+$. This holds by \cite[Corollary 3.2]{pedersen2}.

Now suppose that we have a net $(\tau_i)$ of traces in $\T_I(A)$ converging pointwise on $\ped(I)^+$
to a trace $\tau$, also in $\T_I(A)$. Let $a\in A^+$ and $\epsilon>0$. We need to show that
the inequalities \eqref{topology} hold. If $a\notin I^+$ then this is 
true, since $\tau_i(a)=\tau(a)=\infty$ for all $i$. Suppose that  $a\in I^+$. Then $(a-\epsilon)_+\in \ped(I)^+$
for all $\epsilon>0$. 
So 
\[\limsup \tau_i((a-\epsilon)_+) =\tau((a-\epsilon)_+) \leq \tau(a),\]
and
\[
\tau((a-\epsilon)_+) =\liminf \tau_i((a-\epsilon)_+) \leq \liminf \tau_i(a),
\]
for all $\epsilon>0$. 

(ii) The traces that are a multiple of 0 form a closed subset of $\T(A)$. 
Hence  $\lat(A)$  is compact and Hausdorff
in the relative topology inherited from $\T(A)$. Let us show that this
topology is finer than the 
Fell-Lawson topology. This will give the desired result, since $\lat(A)$
is compact and Hausdorff in both topologies.

Recall that the Fell-Lawson topology has the sub-basis of open sets 
$U_I=\{J\in \lat(A)\mid I\nleq J\}$, and $V_I=\{J\in \lat(A)\mid I\wayb J\}$, 
where $I$ ranges in $\lat(A)$. Here we have denoted by $\wayb$ 
the (countable)  far below relation in the ordered set $\lat(A)$;
see Section 4.2 below (cf.~also \cite{compendium}, where uncountable increasing nets are allowed).
Suppose that $(J_i)_{i\in \Lambda}$ is 
a net converging to $J$ in the relative topology, and $J\in U_I$. If we have $I\leq J_i$ for a subnet $(J_i)_{i\in \Lambda'}$, then $\tau_{J_i}\leq \tau_I$ for all $i\in \Lambda'$, whence $\tau_J\leq \tau_I$. This contradicts the relation 
$I\nleq J$. Therefore, there exists $i_0$ such that $I\nleq J_i$ for all $i\geq i_0$. This shows that the set  $U_I$ is open in the relative topology.

Let $J\in V_I$. For every $b\in J^+$ and $\epsilon>0$ consider the ideal $J_{b,\epsilon}=\Ideal((b-\epsilon)_+)$.
By Proposition \ref{CP} (i), the ideals $J_{b,\epsilon}$ form an upward directed subset of $\lat(A)$. Since they have supremum $J$,
we must have $I\subseteq \Ideal((b-\epsilon)_+)$ for some $b\in J^+$ and $\epsilon>0$. If $(J_i)$ is a net such that $J_i\to J$ in the relative topology,
then $\limsup \tau_{J_i}((b-\epsilon/2)_+)\leq \tau_J(b)=0$. Thus $(b-\epsilon/2)_+\in J_i$ for all $i\geq i_0$, for some $i_0$. 
In other words, $I\subseteq \Ideal((b-\epsilon)_+)\wayb  \Ideal((b-\epsilon/2)_+)\subseteq J_i$ for all $i\geq i_0$. 
This shows that the set $V_I$ is open in the relative topology.
\end{proof}

\subsection{The functor $\T(\cdot)$}
Homomorphisms between C*-algebras induce morphisms in the opposite
direction between their cones of traces; given 
$\phi\colon A\to B$ the map $\T(\phi)\colon \T(B)\to \T(A)$ is defined by $\T(\phi)(\tau)=\tau\circ \phi$. 
It is easily verified that $\T(\phi)$ is linear and continuous.

Let us denote by $\mcC$ the category of compact Hausdorff non-cancellative cones with 
jointly continuous addition and jointly continuous scalar multiplication by $[0,\infty]$, with, as morphisms,
continuous linear maps between cones. 
(Here, linear means additive, homogeneous with respect
to scalars in $[0, \infty]$, and takes $0$ into $0$.)
By Proposition \ref{addition} and 
Theorem \ref{compact-hausdorff}, the cone $\T(A)$ is in the
category $\mcC$. 

\begin{theorem}\label{functor}
$\T(\cdot)$ is a continuous contravariant functor from the category of 
C*-algebras to the category $\mcC$.
\end{theorem}

\begin{proof}
It is straightforward that $\T(\cdot)$ is a functor. (If 
$\phi:A\to B)$ is a homomorphism of C*-algebras and
$\tau_i\to \tau$ in $\T(B)$, then to show that
$\tau_i \phi \to \tau \phi$ in $\T(A)$, passing
to a subnet with $\tau_i \phi$ converging pointwise to $\sigma$ 
we must  show that $\widetilde\sigma=\tau \phi$. We may suppose that
$\tau_i\to \sigma'$ pointwise, so $\sigma=\sigma' \phi$,
and since $(\sigma')^\sim=\tau$, we have
$(\sigma' \phi)^\sim =(\sigma')^\sim \phi=\tau \phi$,
as desired.)
Let $A=\varinjlim (A_i,\phi_{i,j})$ be an inductive limit of C*-algebras. 
N.B: we will not assume that the index set in this inductive
system is countable. Let $C$ denote the subset of the Cartesian product 
$\prod_i \T(A_i)$ of vectors $(\tau_i)$ compatible with the projective 
system $(\T(A_i),\T(\phi_{i,j}))$; that is,
$\tau_i=\T(\phi_{i,j})(\tau_j)$ for all $i<j$. Denote by
$\mu_i\colon C\to \T(A_i)$ the projection
onto the $i$th coordinate. It is well known that $C$ is the projective
limit of $(\T(A_i),\T(\phi_{i,j}))$
in the category of compact Hausdorff spaces. It is easily verified that
$C$ is a cone when endowed
with the operations of coordinate-wise addition and scalar multiplication,
that $C$ belongs to the
category $\mcC$, and that $(C,\mu_i)$
is in fact the projective limit of the system of cones
$(\T(A_i),\T(\phi_{i,j}))$ in the category $\mcC$. 

Let $m\colon \T(A)\to C$ denote the map given by $m(\tau):=(\T(\phi_{i,\infty}(\tau)))$.
In order to show that $\T(A)$ and $C$ are isomorphic,
it is enough to prove that $m$
is bijective, since a continuous bijection between compact Hausdorff spaces has
continuous inverse. 

Suppose that $\tau_1$ and $\tau_2$ are  traces in $\T(A)$ such that $m(\tau_1)=m(\tau_2)$, i.e.,  
$\tau_1\circ\phi_{i,\infty}=\tau_2\circ\phi_{i,\infty}$ for all $i$. Then $\tau_1$ and $\tau_2$ agree
on the set $\bigcup_i \phi_{i,\infty}(A_i^+)$ of positive elements coming from the algebras $A_i$. Let us call this set $B$. We have that $B$ is dense  
in $A^+$ and is such that if $a\in B$ then $(a-\epsilon)_+\in B$ for all $\epsilon>0$. It follows, by  
Remark \ref{afterdef}, that $\tau_1$ and $\tau_2$ cannot be separated
in the topology of $\T(A)$. Since $\T(A)$ is Hausdorff (by Theorem \ref{compact-hausdorff}), this shows that $\tau_1=\tau_2$.

Let $(\tau_i)$ be a vector in $C$. Let us find a trace $\tau$ such that $m(\tau)=(\tau_i)$. 
For $a\in A_i^+$ write $\tau_i(a)=\tau(\phi_{i,\infty}(a))$. 
Let us show that $\tau$ is well defined on the set $B$. Suppose that $\phi_{i,\infty}(a)=\phi_{i,\infty}(b)$. For
every $\epsilon>0$ there exists $j$ such that $\|\phi_{i,j}(a)-\phi_{i,j}(b)\|<\epsilon$.   
Using Lemma \ref{K-R} we get 
\[\tau_i((a-\epsilon)_+)=\tau_j(\phi_{i,j}((a-\epsilon)_+))\leq \tau_j(\phi_{i,j}(b))=\tau_i(b).\]
In the limit as $\epsilon\to 0$  we obtain $\tau_i(a)\leq \tau_i(b)$. Hence by symmetry, $\tau_i(a)=\tau_i(b)$.
So $\tau$ is well defined on $B$.

Let us extend $\tau$ from $B$ to $A^+$ as follows. Define $\widetilde\tau\colon A^+\to [0,\infty]$ by
\[
\widetilde\tau(a):=\sup\{\,\tau(a')\mid a'\in B, a'\preccurlyeq (a-\epsilon)_+ \hbox{ for some }\epsilon>0\,\}.
\]
Let us show that $\widetilde\tau$ is a trace; clearly, it extends $\tau$.
For every $a'\in B$ and $\epsilon>0$ such that $a'\preccurlyeq (a-\epsilon)_+$ we have 
$\widetilde\tau(a')\leq \widetilde\tau((a-\epsilon/2)_+)$. 
It follows from this that $\widetilde\tau(a)=\sup_{\epsilon>0} \widetilde\tau((a-\epsilon)_+)$.
Since $(xx^*-\epsilon)_+\sim (x^*x-\epsilon)_+$ for all $\epsilon>0$ and $x\in A$, we have $\widetilde\tau(xx^*)=\widetilde\tau(x^*x)$.
Also, it can be shown using \eqref{CP1} (see Proposition \ref{CP} (i)) that 
$\widetilde\tau$ is superadditive, i.e., $\widetilde\tau(a)+\widetilde\tau(b)\leq \widetilde\tau(a+b)$.
In particular $\widetilde\tau$ is increasing.

It remains to show that $\tau$ is subadditive. (Homogeneity is clear.)
Before proceeding with the proof of the subadditivity of $\tau$ let us prove a preliminary fact.
Suppose that $c,c'\in A^+$ are such that $\|c-c'\|<\epsilon$. 
By \cite[Lemma 2.2]{kirchberg-rordam} there exists a contraction $d\in A$ such that
$(c-\epsilon)_+=dc'd^*$. Since $\tau$ is increasing and satisfies the trace
identity, it follows that $\widetilde\tau((c-\epsilon)_+)\leq \widetilde\tau(c')$ whenever $\|c-c'\|<\epsilon$. 
 
Let us show that $\tau$ is subadditive.
Let $a,b\in A^+$ and $\epsilon>0$ and assume, as we may without loss of generality, that $\|a\|, \|b\|\leq 1$. 
Choose $a',b'\in B$ such that $\|a-a'\|<\epsilon/4$ and $\|b-b'\|<\epsilon/4$.
Then
\begin{align*}
\widetilde\tau((a+b-\epsilon)_+) &\leq \widetilde\tau((a'-\frac \epsilon 2)_++(b'-\frac \epsilon 2)_+)=\\
&=\widetilde\tau((a'-\frac \epsilon 2)_+)+ \widetilde\tau((b'-\frac \epsilon 2)_+)\leq \widetilde\tau(a)+\widetilde\tau(b).
\end{align*}
Passing to the supremum on the left side with respect to $\epsilon$, we deduce that $\widetilde\tau$ is subadditive. This shows that 
$\widetilde\tau$ is a trace. That $\widetilde\tau$ is lower semicontinuous follows from Lemma \ref{whatsatrace} and the
equation $\widetilde\tau(a)=\sup_{\epsilon>0} \widetilde\tau((a-\epsilon)_+)$.
\end{proof}

\begin{remark}
In addition to being an object in the category $\mcC$, 
we have seen that the cone $\T(A)$
is a complete lattice and satisfies the identities
\eqref{supidentity} and \eqref{infidentity}.
Inspection of the proofs of Propositions \ref{distributive} and
\ref{infdistributivity}
shows that they also hold on replacing $\T(A)$ by any topological
cone in $\mcC$ that is a complete lattice and  satisfies the identities 
\eqref{supidentity} and \eqref{infidentity}.
Is Proposition \ref{inducedtops} (i) true for such cones too? 
\end{remark}

\section{Functionals on the Cuntz Semigroup} 
\subsection{Quasitraces and functionals}
Let $\Cu(A)$ denote the stabilized Cuntz semigroup of $A$, i.e., the ordered semigroup of Cuntz equivalence classes of positive 
elements in  $A\otimes \mathcal K$ (see \cite{coward-elliott-ivanescu} and \cite{rordam}). 
Given $\tau\in \T(A\otimes \mathcal K)$ let us define the function $\lambda_\tau([a]):=\sup_n \tau(a^{1/n})$, where  $[a]\in \Cu(A)$ denotes the equivalence class of $a\in (A\otimes \mathcal K)^+$. This is known to be
a well defined function on $\Cu(A)$ with values in $[0,\infty]$ and with the following properties:

(1) $\lambda_\tau$ is additive and order preserving, and $\lambda_\tau (0)=0$,

(2) $\lambda_\tau$ preserves the suprema of increasing sequences.

Let us refer to a function $\lambda\colon \Cu(A)\to [0,\infty]$ with the properties (1) and (2) 
as a functional on $\Cu(A)$.

\begin{lemma} \label{whatsaquasitrace}
If $\tau:A^+\to [0, \infty]$ is a quasitrace then 
$\widetilde \tau$ defined by $\widetilde\tau(a)=\sup_{\epsilon>0}
\tau((a-\epsilon)_+)$ is a lower semicontinuous quasitrace
on $A$, and is the largest such quasitrace majorized by
$\tau$.
\end{lemma}

\begin{proof}
If $B$ is a commutative sub-C*-algebra of $A$, 
then the restriction $\tau|B$ of $\tau$
to $B$ is a trace on $B$, and the  restriction of $\widetilde \tau$ 
to $B$ is the lower semicontinuous regularization
$(\tau|B)^\sim$ of Lemma \ref{whatsatrace}. 
This shows that $\widetilde\tau$  is additive
on elements that commute.
Since for every $x\in A$ and $\epsilon>0$ there is 
$y\in A$ such that $(x^*x -\epsilon)_+= y^*y$ and
$(xx^*-\epsilon)_+=yy^*$ (see the proof of Proposition \ref{CP}
(ii)), we have $\widetilde\tau(x^*x)=\widetilde\tau(xx^*)$
for all $x\in A$. So $\widetilde \tau$ is a quasitrace. This,
together with the defining equation $\widetilde \tau(a)=\sup_{\epsilon>0}
\tau((a-\epsilon)_+)$, implies that $\widetilde\tau$ is lower
semicontinuous (see the last remark
in \cite[Definition 2.2]{kirchberg-blanchard}). If $\sigma$ is another
lower semicontinuous quasitrace with $\sigma \le\tau$, then (as in the
proof of Lemma \ref{whatsatrace}), for any $a\in A^+$, 
\[
\sigma(a)=\sup_{\epsilon>0}
\sigma((a-\epsilon)_+)\le \sup_{\epsilon>0}\tau(a-\epsilon)_+)=\tau(a).
\qedhere
\]
\end{proof}

To repeat,  we shall denote by $\QT_2 (A)$ the cone of lower
semicontinuous quasitraces of $A\otimes\mathcal K$.
The notation $\QT_2(A)$ is explained by the  result of Blanchard and 
Kirchberg that
every lower semicontinuous 2-quasitrace of $A$ extends to a
lower semicontinuous quasitrace of $A\otimes \mathcal K$;
see \cite[Remark 2.27 (viii)]{kirchberg-blanchard}). Let $\F(\Cu(A))$ denote the
cone of functionals on $\Cu(A)$, as defined above.

\begin{proposition} \label{quasifunct}
Given $\tau\in \QT_2 (A)$ the function $\lambda_\tau([a]):=
\sup_n \tau(a^{1/n})$ is well defined and gives a 
functional on $\Cu(A)$, i.e., an element of $\F(\Cu (A))$.
Given $\lambda \in \F(\Cu(A))$ the function
$\tau_\lambda (a)=\int^\infty_0 \lambda([(a-t)_+])\, dt$ is
a lower semi-continuous quasitrace on $A\otimes \mathcal K$, i.e.,
an element of $\QT_2(A)$. The maps $\tau\mapsto\lambda_\tau$
and $\lambda\mapsto \tau_\lambda$ are the inverses of each other.
\end{proposition}

\begin{proof}
Let $\tau\in \QT_2 (A)$ and set $\sup_{n}\tau (a^{1/n})=D(a)$,
for $a\in (A\otimes \mathcal K)^+$. The
restriction of $\tau$ to a commutative sub-C*-algebra 
$B$ of $A\otimes\mathcal K$ consists of
integration by some measure on the spectrum of $B$. If 
$a\in B^+$ then $D(a)$ is the measure of the open
subset of the spectrum of $B$ consisting of points where
$a$ does not vanish. From this observation it follows that
(1) $D(a+b)=D(a)+D(b)$ if $a$ and $b$ are orthogonal, 
(2) $D(a)\le D(b)$ if $a$ and $b$ commute and 
$a\le Mb$ for some $M>0$, and (3) $D(a)=\sup_{\epsilon>0}
D((a-\epsilon)_+)$.

Let us show that $D$ is constant on the Cuntz equivalence classes
of positive elements. First, note that for any
$x\in A\otimes \mathcal K$ and $g\in C_0 (\R^+)$ there is $y$
such that $g(x^* x)=y^* y$ and $g(xx^*)=yy^*$. Therefore,
$D(x^* x)=D(xx^*)$. Suppose
now that $a$ is Cuntz smaller than $b$, i.e., 
$d^*_n bd_n\to a$ for some sequence $(d_n)$.
Let $\epsilon >0$. Choose $d\in A\otimes \mathcal K$
and $\delta>0$ such that $\|a-d^*(b-\delta)_+d\| <\epsilon$.
By the proof of
Lemma \ref{K-R}, there is $y$ such that $(a-\epsilon)_+=y^*y$
and $yy^*\le (b-\delta)_+$.
Choose a continuous function $g$ with
$g(0)=0$ such that $g(b)\le Mb$ for some $M>0$ and
$g(b)(b-\delta)_+=(b-\delta)_+$. We have 
$g(b)yy^*=yy^*g(b)=yy^*$. In particular,
$yy^*$ and $g(b)$ commute and 
$yy^*\le \|y\|^2 g(b)$. Hence $D(yy^*)\le D(g(b))\le D(b)$, by
(2) of the previous paragraph. Thus, $D((a-\epsilon)_+)=
D(yy^*)\le D(b)$. Letting $\epsilon$ tend to  $0$ and applying (3)
of the previous paragraph we obtain $D(a)\le D(b)$.

It follows from the previous discussion that $\lambda$  defined by 
$\lambda([a]):=D(a)$ is well defined on $\Cu(A)$, additive,
order  preserving, takes $0$ into $0$, and satisfies 
$\lambda([a])=\sup_{\epsilon>0} 
\lambda([(a-\epsilon)_+])$. In order to show that $\lambda$
is a functional it remains to show that it preserves the 
suprema of arbitrary increasing sequences.
Let $([a_n])$ be increasing with supremum $[a]$. It is known
that $[(a-\epsilon)_+] \ll [a]$ (see the first paragraph of
the next subsection for the definition of the relation 
$\ll$---and see \cite{coward-elliott-ivanescu} for the statement $[(a-\epsilon)_+]\ll[a]$).
This implies that $[(a-\epsilon)_+]\le [a_n]$ for
some $n$. Thus, $\lambda([a-\epsilon)_+])\le \sup_n \lambda
[a_n]),$ and letting $\epsilon$ go to $0$ we get
$\lambda([a])\le \sup_n \lambda([a_n])$. The reverse inequality
is clearly true, since $\lambda$ is order preserving.

Let us now start with a functional $\lambda$ and let
$\tau_{\lambda}$ be defined as in the statement of the
proposition. If $B$ is commutative sub-C*-algebra of 
$A\otimes \mathcal K$, and $a\in B^+$, then
$\lambda([a])$ depends only on the set of points in the
spectrum of $B$ where $a$ does not vanish. Moreover,
$\lambda$ defines a Borel measure on the spectrum of $B$
in this way. By Fubini's theorem, $\tau_{\lambda}(a)$ is the
integral of $a$ with respect to that measure. Therefore,
$\tau_{\lambda}$ is additive on $B$ and $\tau_{\lambda}(a)=
\sup_{\epsilon > 0}\tau_{\lambda} ((a-\epsilon)_+)$.

For every $x\in A\otimes \mathcal K$ and $\epsilon>0$ there
exists $y$ such that $(x^* x-\epsilon)_+=y^* y$ and
$(xx^*-\epsilon)_+=yy^*$. Hence $\lambda([x^* x- \epsilon)_+])=
\lambda([xx^*-\epsilon)_+])$, and so
$\tau_\lambda (x^* x)=\tau_\lambda (xx^*)$. It follows that
$\tau_\lambda$ is a quasitrace. We also know that
$\tau_\lambda(a)=\sup_{\epsilon>0} \tau_\lambda((a-\epsilon)_+)$.
This implies that $\tau_\lambda$ is lower semicontinuous 
(see \cite[Definition 2.22]{kirchberg-blanchard}).

Finally, we need to show that the maps $\tau\mapsto \tau_\lambda$ 
and $\lambda\mapsto \lambda_\tau$ are inverse
to each other. It is immediate from the definitions of
these two maps that it is enough to prove this on the
commutative sub-C*-algebra generated by a positive element. In
this case the result follows from standard results in the theory 
of integration.
\end{proof}

\begin{remark}\label{haagerup}  
A theorem of Haagerup says that if $A$ is exact and unital then
every bounded 2-quasitrace on $A$ is a trace. It was observed in 
\cite[Remark 2.29 (i)]{kirchberg-blanchard} that after a number of elementary reductions
this theorem can be extended to obtain that every lower
semicontinuous 2-quasitrace on an exact C*-algebra must be a trace.
It follows that if $A$ is exact then every functional on $\Cu(A)$ arises from a lower semicontinuous 
trace. 
\end{remark}

Let us  endow $\QT_2(A)$ with the topology in which the net $(\tau_i)$
converges to $\tau$ if
\begin{equation}\label{QTtopology}
\limsup \tau_i((a-\epsilon)_+)\leq \tau(a)\leq \liminf \tau_i(a)
\end{equation}
for all $a\in (A\otimes \mathcal K)^+$ and $\epsilon>0$. Alternatively
(as for $\T(A)$), $\tau_i\to \tau$ if, 
whenever a subnet of $(\tau_i)$ converges pointwise to a function
$\sigma\colon (A\otimes \mathcal K)^+\to [0,\infty]$,
the regularization $\widetilde\sigma$ of $\sigma$ given by Lemma \ref{whatsaquasitrace}
above is equal to $\tau$.
A neighbourhood basis of a point of $\QT_2(A)$ 
(or for a uniform structure determining the topology)
can be described just as
before for $\T(A)$. Notice that the relative topology 
of $\T(A\otimes \mathcal K)$, as a subset of $\QT_2(A)$, is the same
topology that we assigned to
$\T(A\otimes \mathcal K)$ in the previous section. 

Let us say that
the net $(\lambda_i)$ converges to  $\lambda$ on the cone
$\F(\Cu(A))$ of functionals on $\Cu(A)$  if
\[
\limsup \lambda_i([(a-\epsilon)_+])\leq \lambda([a])\leq
\liminf \lambda_i([a])
\] 
for all $[a]\in \Cu(A)$ and $\epsilon>0$. Alternatively
(as for $\T(A)$), $\lambda_i\to\lambda$
if $\lambda([a])=\sup_{\epsilon>0}\lambda'([(a-\epsilon)_+])$ for all
$a\in (A\otimes \mathcal K)^+$, whenever
a subnet of $(\lambda_i)$ converges pointwise to a function
$\lambda'\colon \Cu(A)\to [0,\infty]$. 

\begin{theorem}\label{isocompacts}
The cones $\QT_2(A)$ and $\F(\Cu(A))$ are compact and Hausdorff, and
the map $\tau\mapsto \lambda_\tau$ is a a homeomorphism between them.
\end{theorem}

\begin{proof}
The proof that $\QT_2(A)$ is compact and Hausdorff is similar to the  proof
given above for $\T(A)$
(Theorem \ref{compact-hausdorff}). 
This is also the case for the proof  that $\F(\Cu(A))$ is compact and Hausdorff 
(see Theorem \ref{Dcompact} below for a generalization of this). 
In order to show that $\tau\mapsto \lambda_\tau$ is a homeomorphism it is
enough to show that it is continuous.
Let $(\tau_i)$ be a net in $\QT_2(A)$ such that $\tau_i\to \tau$ in $\QT_2(A)$. 
Let $[a]\in \Cu(A)$.  We have 
\[
\tau(a^{1/n})\leq \liminf \tau_i(a^{1/n})\leq \liminf \lambda_{\tau_i}([a]),
\] 
for all $n$. Therefore, $\lambda_\tau([a])\leq \liminf \lambda_{\tau_i}([a])$.

Let $\epsilon>0$.
Choose $f\in C_0(\R^+)$ such that $0\leq f\leq 1$
and $f(t)=1$ for $t\in [\epsilon,\|a\|]$. Set $f(a)=a'$.
We have 
\[
(1-\epsilon)\limsup \lambda_{\tau_i}([(a-\epsilon)_+])
\leq \limsup \tau_i(a'-\epsilon)_+ \leq \tau(a')\leq \lambda_\tau([a]),
\]
for all $\epsilon>0$.  This implies that 
$\limsup \lambda_{\tau_i}([(a-\epsilon)_+])\leq \lambda_\tau([a])$.
\end{proof}

\begin{remark} If $A$ is exact then $\T(A\otimes \mathcal K)=\QT_2(A)$
(see Remark \ref{haagerup}).
Since it is always the case that $\T(A)\cong \T(A\otimes \mathcal K)$, we have that 
if $A$ is exact then $\T(A)\cong \T(A\otimes \mathcal K)=\QT_2(A)\cong \F(\Cu(A))$. 
\end{remark}

\subsection{The category $\mcCu$}
In \cite{coward-elliott-ivanescu}, 
Coward, Elliott, and Ivanescu showed that $\Cu(A)$ belongs to a 
particular category of ordered semigroups denoted by $\mcCu$. Let us 
recall the definition of this category here.

For elements $a$ and $b$ of an ordered set, let us say that $a$ is 
far below $b$, and write 
$a\wayb b$, if for any increasing sequence $(b_n)$ with supremum
greater than or equal to $b$ there 
exists $n$ such that $a\leq b_n$. (Then in particular $a\leq b$.)
The category $\mcCu$ has 
for objects the ordered semigroups $S$ with 0 such that

(1) increasing sequences in $S$ have a supremum,

(2) for every  $a\in S$ there is a sequence $a_1,a_2,\dots$ with supremum $a$
such that $a_n\wayb a_{n+1}$ for all $n$,

(3) if $a_1\wayb b_1$ and $a_2\wayb b_2$ then $a_1+a_2\wayb b_1+b_2$, and

(4) if $(a_n)$ and $(b_n)$ are increasing sequences then
$\sup(a_n+b_n)=\sup a_n+\sup b_n$.

The morphisms of the category $\mcCu$ are the ordered semigroup morphisms
(i.e., the additive and order 
preserving maps) that preserve suprema of increasing sequences and the far
below relation. 

\begin{remark} The far below relation (also referred to as the way below relation, or, more formally, compact containment), is usually defined with respect to increasing nets
$(b_i)$ instead of increasing (countable) sequences $(b_n)$. Nevertheless, it is (countable) increasing sequences that we wish to consider here. To avoid confusion we might say countable compact containment. 
\end{remark}

Let $S$ be a semigroup in the category
$\mcCu$. Let us call functionals on $S$ those additive and order preserving  
functions from $S$ to $[0,\infty]$ that take 0 into 0 and preserve the suprema of
increasing sequences.
Let us denote by $\F(S)$ the cone of  functionals on $S$ endowed with pointwise
addition and scalar multiplication
by strictly positive real numbers. Let us consider $\F(S)$ with the topology in 
which  a net $(\lambda_i)$ converges to 
a point $\lambda$ if 
\begin{align}\label{Dtopology}
\limsup \lambda_i(x)\leq \lambda(y)\leq \liminf \lambda_i(y) 
\end{align}
for all $x,y\in S$ such that $x\wayb y$. If $\phi\colon S\to T$ is a
morphism in the category 
$\mcCu$ then $\F(\phi)(\lambda):=\lambda\circ \phi$, $\lambda\in \F(T)$, is
a continuous linear
map from $\F(T)$ to $\F(S)$.  

\begin{lemma}\label{whatsaD}
Let $S$ be a semigroup in the category $\mcCu$ and let 
$\lambda\colon S\to [0,\infty]$ 
be an additive map on $S$. Then $\widetilde\lambda$ defined by
$\widetilde \lambda(x):=\sup\{\,\lambda (z)\mid z\wayb x\,\}$
is a functional on $S$, and is the largest functional majorized by
$\lambda$.
\end{lemma}
\begin{proof}
If $x\leq y$ and $x'\wayb x$ then $x'\wayb y$. This allows us to conclude
that $\widetilde \lambda(x)\leq \widetilde\lambda(y)$.
If $x'\wayb x$ and $y'\wayb y$ then $x'+y'\wayb x+y$, and so 
$\widetilde\lambda(x)+\widetilde\lambda(y)\leq \widetilde\lambda(x+y)$.
If $z\wayb x+y$ then there are $x'\wayb x$ and $y'\wayb y$ such that
$z\wayb x'+y'$. It follows that
 $\widetilde\lambda(x+y)\leq \widetilde\lambda(x)+\widetilde\lambda(y)$.
Finally, let $(x_n)$ be an increasing sequence with supremum $x$. Since
$\widetilde\lambda$
is order preserving we have $\sup \widetilde\lambda(x_n)\leq \widetilde\lambda(x)$.
On the other hand, for all $x'\wayb x$ we have $x'\wayb x_n\leq x$ for some $n$.
Therefore, 
$\widetilde\lambda(x)\leq \sup \widetilde \lambda(x_n)$.  
The last statement is proved as in \ref{whatsatrace} and \ref{whatsaquasitrace}.
\end{proof}

The order of pointwise comparison of functionals in $\F(S)$
is the same as the order arising from the semigroup structure. The proof of
this is identical to the proof for 
$\T(A)$ (see Proposition \ref{orderoftraces}), provided that Lemma \ref{whatsaD}
is used instead of Lemma \ref{whatsatrace}.  

Just as with $\T(A)$, we may extend the scalar multiplication of $\F(S)$
to include 0 and $\infty$:
\begin{align}
(\infty \cdot \lambda)(x)&:=0\hbox{ if }\lambda(x)=0,
\hbox{ and }(\infty \cdot \lambda)(x):=\infty\hbox{ otherwise},\\ 
(0\cdot\lambda)(x)&:=0\hbox{ if }\lambda(z)<\infty, \forall z\wayb x,
\hbox{ and }(0\cdot \lambda)(x)=\infty\hbox{ otherwise.} \label{multby0}
\end{align}
Notice that $0\cdot \lambda$ is the regularization (as in the statement of
Lemma \ref{whatsaD})
of the additive map $x\mapsto 0\cdot \lambda(x)$, where $0\cdot \infty$ is
taken to be $\infty$.
Notice also that $\alpha\lambda\to 0\cdot\lambda$ and
$\frac 1\alpha\lambda\to\infty\cdot \lambda$ as 
$\alpha\to 0$, and, indeed, extended scalar multiplication is (jointly)
continuous overall.

\begin{theorem}\label{Dcompact} 
$\F(\cdot)$ is a sequentially continuous contravariant functor from the
category $\mcCu$ to the category of topological cones $\mcC$.
\end{theorem}

\begin{proof}
The proof that $\F(S)$ is compact and Hausdorff is similar to the proof for
$\T(A)$ (cf.~also Theorem \ref{isocompacts}). 
We use Lemma \ref{whatsaD} instead of Lemma \ref{whatsatrace}. 

Let us show that $\F(\cdot)$ is a sequentially continuous functor. 
As can be seen from the construction given in \cite{coward-elliott-ivanescu},
inductive limits in the category $\mcCu$ are characterized as follows: $S$ is
the inductive limit of 
$(S_i,\phi_{i,j})$ if

(1) every element of $S$ is supremum of an increasing sequence of elements coming
from the $S_i$s,

(2) if $x,y\in S_i$ are such that $\phi_{i,\infty}(x)\leq \phi_{i,\infty}(y)$ in
the limit, 
then for all $z\wayb x$ in $S_i$ there is $n\geq i$ such that
$\phi_{i,n}(z)\leq \phi_{i,n}(y)$ in $S_n$.

Let $C$ denote the projective limit of $(\F(S_i),\F(\phi_{i,j}))$
in the category $\mcC$ (cf. proof of Theorem \ref{functor}).
If $\lambda_1$ and $\lambda_2$ are two functionals on $S$ that agree on
the elements coming
from finite stages, then  $\lambda_1$ and $\lambda_2$ are equal by
the property (1) above 
of inductive limits in $\mcCu$. Thus, the map from $C$ to $\F(S)$ is in injective.
In order to see
that this map is surjective we need to show that for any sequence of
functionals $\lambda_i\in \F(S_i)$
compatible with the inductive limit, there is $\lambda\in \F(S)$ such
that $\F(\phi_{i,\infty})(\lambda)=\lambda_i$.
Let us define $\lambda$ on the subsemigroup $\bigcup_i\phi_{i,\infty}(S_i)$
of $S$ by $\lambda(\phi_{i,\infty}(x))=\lambda_i(x)$.
Let us check that this map is well defined. Suppose that
$\phi_{i,\infty}(x)=\phi_{i,\infty}(y)$. Then by the property (2)
of  inductive limits in the category  $\mcCu$, for every $z\wayb x$ there is $n$ such that 
$\phi_{i,n}(x)\leq \phi_{i,n}(y)$. So
$\lambda_i(z)=\lambda_n(\phi_{i,n}(z))\leq \lambda_n(\phi_{i,n}(y))=\lambda_i(y)$.
Since this holds for all $z\wayb x$ we have $\lambda_i(x)\leq \lambda_i(y)$, 
whence, by symmetry,
$\lambda_i(x)= \lambda_i(y)$.

Let us write $T=\bigcup_i \phi_{i,\infty}(S_i)$. Let us extend $\lambda$
from $T$ to all of $S$ as follows:
\[
\widetilde\lambda(x)=\sup \{\,\lambda(x')\mid x'\wayb x, x'\in T\,\}.
\]
One can now show that $\widetilde\lambda$ is a functional on $S$ that
extends $\lambda$. We will only show here that $\widetilde\lambda$
is additive. Let $x,y\in S$. Let $x'\wayb x$, $y'\wayb y$ and $x',y'\in T$.
Then $x'+y'\wayb x+y$ and $x'+y'\in T$. This 
implies that $\widetilde\lambda$ is superadditive. On the other hand,
if $z'\wayb x+y$ , $z'\in T$, then there are $x',y'\in T$ such that
$z\leq  x'+y'\wayb x+y$ and $x'\wayb x$, $y'\wayb y$. From this we
conclude that $\widetilde\lambda$ is subadditive.  
 \end{proof}

\begin{remark} It was shown in \cite{coward-elliott-ivanescu} that 
$\Cu(\cdot)$
is a sequentially continuous covariant functor
from the category of C*-algebras to the category $\mcCu$. Therefore, by
Theorem \ref{Dcompact}, 
$\F(\Cu(\cdot))$ is a sequentially continuous contravariant functor from
the category of C*-algebras to the
category $\mcC$.
\end{remark}

\section{Dual cones for  $\F(\Cu(A))$ and $\T(A)$}\label{duals}
\subsection{The space $\LL(\F(\Cu(A)))$} Before discussing the dual cone of 
the cone $\F(\Cu(A))$ 
let us begin with some general considerations concerning the cones in
the category $\mcC$.

Let $C$ be a cone in the category $\mcC$. 
Let $\LS(C)$ denote the set of lower semicontinuous functions on $C$
with values
in $[0,\infty]$ that are additive, homogeneous (with respect to the
scalar multiplication by $\R^+$), and  take 0 into 0. We shall
regard $\LS(C)$ as a non-cancellative cone endowed with the operations
of pointwise addition 
and pointwise multiplication by strictly positive scalars. We shall
also consider $\LS(C)$ as ordered by the order of pointwise comparison
of functions. Notice that $\LS(C)$ is closed under passage to suprema
of upward directed sets.
Thus, we can extend the scalar multiplication to include $\infty$ by 
setting $\sup_n n\cdot f=\infty \cdot f$.
For functions $f$ and $g$ in $\LS(C)$ we shall write $f\wayb g$ if for every
increasing sequence $(g_n)$ such that $g\leq \sup g_n$ there is $g_{n_0}$
such that $f\leq g_{n_0}$. 

Let us denote by $\LL(C)$ the subset of $\LS(C)$ composed of those
functions $f$ for which
there is an increasing sequence $(h_n)$, $h_n\in \LS(C)$, with the
following two properties:

(I) the supremum of the $h_n$s is $f$,

(II) $h_n$ is continuous at each point where $h_{n+1}$ is finite.

The definition of $\LL(C)$ is motivated by Proposition \ref{a_is_sta} 
below. It is shown there
that the functions arising from the positive elements of a C*-algebra $A$
(and of $A\otimes \mathcal K$) on the cones $\T(A)$ and
$\F(\Cu(A))$ satisfy (I) and (II) (for suitable increasing sequences).

For $f\in \LS(C)$ let us write $\Set(f)=\{\,\lambda\in C\mid f(\lambda)>1\,\}$.
Notice that
$f\leq g$ if and only if $\Set(f)\subseteq \Set(g)$.

\begin{proposition}\label{closedbysups} Let $f$ and $g$ be in $\LS(C)$. 
 
(i) If $f\leq (1-\mu)g$ for some $\mu>0$, and $f$ is continuous at each 
point where $g$ is finite, then 
$\overline{\Set(f)}\subseteq \Set(g)$. 

(ii) If $\overline{\Set(f)}\subseteq \Set(g)$ then $f\wayb g$ in $\LS(C)$.

(iii) $\LL(C)$ is a subcone of $\LS(\F(\Cu(A)))$ closed under passage to
suprema of increasing sequences.
\end{proposition}

\begin{proof}
(i) Let $(\lambda_i)$ be a net in $C$ such that
$\lambda_i\to \lambda$ and $\lambda_i\in \Set(f)$. 
If $g(\lambda)<\infty$ then $f(\lambda)=\lim f(\lambda_i)\geq 1$, and
so $g(\lambda)>1$. 

(ii) Let  $(g_n)$ be an increasing sequence of functions in $\LS(C)$ with
pointwise supremum
greater than or equal to $g$. Then $\Set(g)\subseteq \bigcup_n \Set(g_n)$.
Since $\overline{\Set(f)}$ 
is compact, we must have that $\overline{\Set(f)}\subseteq \Set(g_{n_0})$
for some $n_0$. Therefore,
$f\leq g_{n_0}$. 

(iii) It follows easily from its definition that $\LL(C)$ is closed under
addition and 
multiplication by strictly positive scalars.

Let $(f_n)$ be an increasing sequence of functions in $\LL(C)$ with
supremum $f$. For every $f_n$ let
$(h_k^n)_{k=1}^\infty$ be a choice of the corresponding sequence
satisfying (I) and (II). 
We may assume without loss of generality that $h_k^n\leq (1-\mu_k^n)h_{k+1}^n$
for some $\mu_k^n>0$.
Let $k_1=1$. Since $h_{k_1+1}^1\wayb f_2$ on $\LS(C)$ (by (i) and (ii)),
there is $h_{k_2}^2$
such that  $h_{k_1+1}^1\leq h_{k_2}^2$. In the same way we may find
$h_{k_3}^3$ such that  $h_{k_2+1}^2\leq h_{k_3}^3$. 
Continue in this way to obtain a sequence $(h_{k_n}^n)$ such that
$h_{k_n+1}^i\leq h_{k_{n+1}}^{n+1}$. 
By proceeding as in the proof 
of \cite[Theorem 1 (i)]{coward-elliott-ivanescu}, we can choose
this sequence so that its  supremum is $f$. If 
$h_{k_{n+1}}^{n+1}(\lambda)<\infty$ then $h_{k_n+1}^n(\lambda)<\infty$, 
and it follows that $h_{k_n}^n$
is continuous at $\lambda$. This shows that $f$ belongs to $\LL(C)$.
\end{proof}

\begin{remark}
The cone $\LL(C)$  is not to be confused with a cone belonging to
the category $\mcC$. In particular,
no topology will be defined on $\LL(C)$. Instead, we shall consider
$\LL(C)$ as a non-cancellative
cone endowed with an order---that of pointwise comparison of
functions---which may not coincide with
the order arising from the addition operation of $\LL(C)$. Also,
we shall not define a scalar multiplication
by $0$ in $\LL(C)$.
\end{remark}

Let us now specialize the study of $\LL(C)$ to the case that
$C=\F(\Cu(A))$ for some C*-algebra $A$.
Our main result is Theorem \ref{supsofas}. As applications of
this theorem we will obtain that $\LL(\F(\Cu(A)))$ is an
ordered semigroup in the category $\mcCu$ and  that $\LL(\F(\Cu(\cdot)))$
is a sequentially continuous functor from the category of C*-algebras to
the category $\mcCu$. 
We will also make use of this theorem in the next section when we
look at the structure 
of the Cuntz semigroup for certain C*-algebras.

Let $A$ be a C*-algebra and $a$  be a positive element of
$A\otimes \mathcal K$. The Cuntz
semigroup element $[a]$ and the positive element $a$ give rise to functions
on $\F(\Cu(A))$:
\begin{align}
\widehat{[a]}(\lambda) &:=\lambda([a]),\\
\widehat{a}(\lambda) &:= \tau_{\lambda}(a)=\int_0^\infty \lambda([(a-t)_+])\, dt,\label{ahat}
\end{align}
where $\lambda\in \F(\Cu(A))$ and $\tau_\lambda$ is the quasitrace associated
to $\lambda$ by Proposition \ref{quasifunct}.
The function $\widehat{[a]}$ belongs to $\LS(\F(\Cu(A)))$ by 
the inequalities \eqref{Dtopology} that define the topology on $\F(\Cu(A))$).
The function
$\widehat{a}$ belongs to $\LS(\F(\Cu(A)))$ by the inequalities that define
the topology on $\QT_2(A)$
and the isomorphism between $\F(\Cu(A))$ and $\QT_2(A)$.

For the rest of this section if $a$ is a positive element of a C*-algebra
we will use the notation 
$a_\epsilon$ to mean the positive element $(a-\epsilon)_+$.

\begin{proposition}\label{a_is_sta}
For all $a\in A^+$ and $\epsilon>0$ the function 
$\widehat{a_\epsilon}$ is continuous at each point where $\widehat a$ is finite.
\end{proposition}

\begin{proof}
Let $\lambda$ be such that $\widehat a(\lambda)<\infty$ and let  
$(\lambda_i)$ be a net in $\F(\Cu(A))$ that converges to $\lambda$. 
We have $\tau_\lambda(a_\epsilon)\leq \liminf
\tau_{\lambda_i}(a_\epsilon)$.
Let $\mu>0$ and set $a'=a_\epsilon + \mu a$. There is $\epsilon'>0$
such that 
$a_\epsilon\leq a_{\epsilon'}'$ (this is easily verified in $C^*(a)$).
Therefore, 
$\limsup \tau_{\lambda_i} (a_\epsilon)\leq \tau_\lambda(a')=
\tau_\lambda(a_\epsilon)+ \mu \tau_\lambda(a)$. This is true for all $\mu>0$.
Since $\tau_\lambda(a)$ is finite we conclude that 
$\limsup \tau_{\lambda_i}(a_\epsilon)\leq \tau_\lambda(a_\epsilon)$.
\end{proof}

Proposition \ref{a_is_sta} implies that $\widehat a$ is in $\LL(\F(\Cu(A)))$
for every $a\in (A\otimes \mathcal K)^+$. 
Since $\widehat{[a]}=\sup_n \widehat{(a^{1/n})}$ and $\LL(\F(\Cu(A)))$ is
closed under passage to suprema of increasing sequences, 
$\widehat{[a]}$ is in $\LL(\F(\Cu(A)))$ too. Notice also that since 
$\widehat{a_\epsilon}\leq (1-\epsilon)\widehat a$, we have that
$\widehat{a_\epsilon}\wayb \widehat a$, by Propositions \ref{closedbysups}
and \ref{a_is_sta}.

Let $I$ be a closed two-sided ideal of $A\otimes \mathcal K$.
Let $f_I\colon \F(\Cu(A))\to [0,\infty]$ denote the
function given by 
\begin{align}
f_I(\lambda)=\left \{
\begin{array}{ll}
0 & \hbox{if }\lambda([a])=0 \hbox{ for all }a\in I^+,\\
\infty & \hbox{otherwise}.
 \end{array}
\right.
 \end{align}
It can be verified that  $f_I$ is in $\LS(\F(\Cu(A)))$. Moreover, every
function in $\LS(\F(\Cu(A)))$ with the only possible
values $0$ and $\infty$ has the form $f_I$ for some ideal $I$.
For $f\in \LS(\F(\Cu(A)))$ let us
write $\Ideal(f)$ for the ideal of $A\otimes \mathcal K$ such that
$\infty\cdot f=f_{\Ideal(f)}$.  If $a$ is a positive element then 
$\Ideal(\widehat a)$ is the closed two-sided ideal generated by
$a$ (i.e., $\Ideal(a)$).

\begin{lemma}\label{wayb} 
Let $f,g\in \LS(\F(\Cu(A)))$. If $\overline{\Set(f)}\subseteq \Set(g)$ then
$\Ideal(f)\wayb \Ideal(g)$ (in $\lat(A)$).
\end{lemma}
\begin{proof}
Let $(I_i)$ be an upward directed collection of ideals with supremum $\Ideal(g)$.
The functions $f_{I_i}$ form an upward directed subset of $\LS(\F(\Cu(A)))$
with supremum $\infty \cdot g$, and so 
$\Set(\infty\cdot g)=\bigcup_i \Set(f_{I_i})$. Since $\overline{\Set(f)}$
is compact and $\overline{\Set(f)}\subseteq
\Set(g)\subseteq \Set(\infty\cdot g)$ (as $g\leq \infty\cdot g$),
we have $\Set(f)\subseteq \Set(f_{I_{i_0}})$ for some $i_0$.
From this we get $\infty \cdot f\leq f_{I_{i_0}}$ that  is to say,
$\Ideal(f)\subseteq I_{i_0}$. 
\end{proof}

The following proposition relies on a result in the duality theory of
topological vector spaces.

\begin{proposition}\label{hahnbanach}
Let $\F_A(\Cu(A))$ denote the subcone of $\F(\Cu(A))$ of functionals 
such that $0\cdot\lambda=0$.
Let $\mathrm{V}(\F_A(\Cu(A)))$ denote the ordered vector space
of linear, real-valued, continuous functions on $\F_A(\Cu(A))$. Then for 
every  positive linear functional 
$\Lambda\colon \mathrm{V}(\F_A(\Cu(A)))\to \R$  there is
$\lambda\in \F_A(\Cu(A))$ such that $\Lambda(f)=f(\lambda)$.
\end{proposition}
\begin{proof}
By \eqref{multby0}, $0\cdot\lambda=0$ if and only if
$\lambda([a_\epsilon])<\infty$ for all $a\in (A\otimes \mathcal K)^+$
and $\epsilon>0$.
Let us identify---via 
Proposition \ref{quasifunct}---the subcone $\F_A(\Cu(A))$ with
the quasitraces of $A\otimes \mathcal K$ 
that are densely finite.  Notice that $\F_A(\Cu(A))$ is a cancellative
cone, since
the quasitraces in it are densely finite.

Let us show that the  relative topology on $\F_A(\Cu(A))$ induced by the 
topology of $\F(\Cu(A))$ is the topology of 
pointwise convergence on the set
$\{\, a_\epsilon\mid a\in (A\otimes \mathcal K)^+, \epsilon>0\, \}$. 
Let $(\lambda_i)$ be a net in $\F_A(\Cu(A))$ and $\lambda\in \F_A(\Cu(A))$. 
Suppose that $\lambda_i\to \lambda$ in the topology of $\F(\Cu(A))$. Since 
$\widehat{a_{\epsilon/2}}(\lambda)<\infty$ for all
$a\in (A\otimes \mathcal K)^+$ and
$\epsilon>0$, by Proposition \ref{a_is_sta}
$\widehat{a_\epsilon}(\lambda_i)\to \widehat{a_\epsilon}(\lambda)$.
Suppose on the other hand that 
$\widehat {a_\epsilon}(\lambda_i)\to \widehat {a_\epsilon}(\lambda)$
for all $a\in A^+$ and $\epsilon>0$. Then we may proceed as in
the proof of  Proposition \ref{inducedtops} (i) to conclude
that $\lambda_i\to \lambda$ in the relative topology of 
$\F_A(\Cu(A))$. 

We conclude that the relative topology on $\F_A(\Cu(A))$ is the
weak*-topology of pointwise convergence on the set 
$\{\,a_\epsilon\mid a\in (A\otimes \mathcal K)^+\,\}$. Therefore,
$\F_A(\Cu(A))$
is a weakly complete cancellative cone in the class $\mathcal S$
of Choquet (see \cite[page~194]{choquet}). 
Now the theorem follows from \cite[Proposition 30.7]{choquet}. 
\end{proof}

\begin{lemma}\label{relunit}
Let $h_1,h_2,h_3\in \LS(\F(\Cu(A)))$ be such that $h_i\leq (1-\mu_i)h_{i+1}$
and $h_i$ is continuous 
at each point where $h_{i+1}$ is finite, for $i=1,2$ and some $\mu_1,\mu_2>0$. 
Then for every $\delta>0$ there is $a\in (A\otimes \mathcal K)^+$
such that $\widehat a\leq h_3$ and $h_1\leq \delta h_3+\widehat a$. 
\end{lemma}
\begin{proof}
We may assume without loss of generality that
$\Ideal(h_3)=A\otimes \mathcal K$, i.e., $h_3(\lambda)>0$
unless $\lambda=0$. Then we have $\widehat b\leq \infty \cdot h_3$ for every
$b\in (A\otimes \mathcal K)^+$.
Since $\widehat{b_{\epsilon}}\wayb \widehat b$ for every $\epsilon>0$,
we obtain that
$\widehat{b_{\epsilon}}\leq Mh_3$ for some finite $M>0$. 
Set $K=\{\,\lambda\in \F(\Cu(A))\mid h_3(\lambda)\leq 1\,\}$. 
By Proposition \ref{a_is_sta}, the function  $\widehat{b_\epsilon}$ 
is continuous on $K$ for all $b\in (A\otimes \mathcal K)^+$ and all $\epsilon>0$.
Also, by hypothesis, $h_1$ is continuous on $K$. 

Let us show that $h_1$ can be uniformly approximated on $K$ by convex combinations
of functions of the form $\widehat {b_\epsilon}$. 
Suppose the contrary. Then there is a real 
measure $m$ on $K$ such that $\int \widehat{b_\epsilon}\, dm=0$ for all
$b\in (A\otimes \mathcal K)^+$ and $\epsilon>0$, and 
$\int h_1\,dm=1$. Let $m=m_+-m_-$ denote the Jordan decomposition of $m$.  Then
\[
\int h_1\,dm_+=\int h_1\,dm_-+1 \hbox{\quad and \quad}
\int \widehat {b_\epsilon} \,dm_+=\int \widehat {b_\epsilon} \,dm_-,\hbox{\quad for all }b,\epsilon.
\] 
For $\lambda\in K$ we have $h_3(0\cdot \lambda)=0$, and since
$\Ideal(h_3)=A\otimes \mathcal K$, we conclude that 
$0\cdot \lambda=0$. Therefore, $K$ is contained in $\F_A(\Cu(A))$, 
with $\F_A(\Cu(A))$ as defined in Proposition \ref{hahnbanach}.
So we can define positive linear functionals $\Lambda_+$, $\Lambda_{-}$ on
the vector space $\mathrm{V}(\F_A(\Cu(A)))$
by $\Lambda_+(f)=\int fd\,m_+$
and $\Lambda_{-}=\int f\,dm_{-}$. By Proposition \ref{hahnbanach},
$\Lambda_+$ and $\Lambda_-$ are given by evaluation  on functionals 
$\lambda_+$ and $\lambda_-$  belonging to $\F_A(\Cu(A))$. 
For all $b\in (A\otimes \mathcal K)^+$ and $\epsilon>0$ we have 
$b_\epsilon\in \mathrm{V}(\F_A(\Cu(A)))$. 
So $\widehat{b_\epsilon}(\lambda_+)=\widehat {b_\epsilon}(\lambda_-)$
for all $b$ and $\epsilon$. This implies that 
$\lambda_{+}=\lambda_{-}$. 

Let us show that the restriction of $h_1$ to $\F_A(\Cu(A))$
is also in $\mathrm{V}(\F_A(\Cu(A)))$. 
Let $\lambda\in \F_A(\Cu(A))$, i.e., $0\cdot\lambda=0$.
If $h_2(\lambda)=\infty$ then
$h_2(\lambda/n)=\infty$ for all $n$. Since 
$\overline{\Set(h_2)}\subseteq \Set(h_3)$ (by Proposition \ref{closedbysups} (i)), 
$0\cdot \lambda\in \Set(h_3)$. This implies that $h_3(0\cdot \lambda)=\infty$,
which  contradicts the equation $0\cdot\lambda=0$. 
We conclude that $h_2(\lambda)<\infty$ for all $\lambda\in \F_A(\Cu(A))$,
and so  $h_1$ is continuous on $\F_A(\Cu(A))$.

We now have $h_1(\lambda_+)=h_1(\lambda_{-})+1$. This  contradicts the
earlier conclusion $\lambda_+=\lambda_-$.
Therefore, the restriction of  $h_1$ to $K$ belongs to the closure of
the convex set spanned by the functions 
$\widehat{b_\epsilon}$. Hence, for
every $\delta>0$ there exists a positive element $a$ such that
$\|h_1-a\|_K<\delta$. Equivalently,
$h_1\leq \widehat a+\delta h_3$
and $\widehat a\leq h_1+\delta h_3$ on $K$. It is easily shown that
these inequalities also  hold on all $\F(\Cu(A))$. 
Changing $\widehat a$ to $\widehat a/(1+\delta)$ we can arrange that
$\widehat a\leq h_3$.
\end{proof}

\begin{theorem} \label{supsofas}
Let $f$ be in $\LL(\F(\Cu(A)))$. Then $f$ is the supremum of an
increasing sequence
of $\widehat a$s ($a\in (A\otimes \mathcal K)^+$). Such a sequence may
even be chosen to be rapidly increasing: $a_1 \wayb a_2 \wayb \cdots$.
\end{theorem}
\begin{proof}
Let $(h_n)$ be a an increasing sequence satisfying (I) and (II). We may
assume without 
loss of generality that $h_n\leq (1-\mu_n)h_{n+1}$ for some $\mu_n>0$,
for all $n$.
By Proposition \ref{closedbysups}, $\overline{\Set(h_n)}\subseteq \Set(h_{n+1})$
and $h_n\wayb h_{n+1}$ for all $n$. 
Hence, by Lemma \ref{wayb}, $\Ideal(h_n)\wayb \Ideal(h_{n+1})$
in $\lat(A)$ for all $n$. Let us choose $b\in (A\otimes \mathcal K)^+$
and $\epsilon_0>0$  such that 
$\Ideal(h_4)\subseteq \Ideal(b_{\epsilon_0})$ and
$\Ideal(b)\subseteq \Ideal(h_5)$. 
We have that $h_4\leq \infty \cdot \widehat {b_{\epsilon_0}}$ and
$\widehat b\leq \infty \cdot h_5$. 
Therefore, there is a  constant $M>0$ such that
$h_2\leq  M \widehat {b_{\epsilon_0}}$ and $\widehat {b_{\epsilon_0}}\leq  M h_5$. 
Let us choose $\delta$ such that $\delta M < \mu_3$. Finally,
using Lemma \ref{relunit}, let us find 
$a$ in $(A\otimes \mathcal K)^+$ such that $\widehat a\leq h_3$ and
$h_1\leq (\delta/M) h_3+\widehat a$. 

By the stability of $A\otimes \mathcal K$, we may assume that the positive
elements $a$ and $b$ that we found in the
previous paragraph are orthogonal to each other. (If they are not, we
may replace them by Murray-von Neumann equivalent
elements that are orthogonal.) Let $a_1=a+\delta b_{\epsilon_0}$. Then 
\[
\widehat a_1= \widehat a+\delta\widehat {b_{\epsilon_0}}\leq
(1-\mu_3+\delta M)h_5\leq h_5.
\]
Also 
\[
\widehat a_1= \widehat a+\delta\widehat {b_{\epsilon_0}}\geq
\widehat a+\frac \delta M h_3\geq h_1.
\] 
So $h_1\leq \widehat a_1\leq h_5$. In the same way we may find
$\widehat a_2$ such that
$h_5\leq \widehat a_2\leq h_9$. Continuing in this way we get the
desired sequence---by Proposition \ref{closedbysups} (i) and (ii) even rapidly increasing.
\end{proof}

\begin{corollary}\label{supsof[a]}
Every function in $\LL(\F(\Cu(A)))$ is the supremum of an increasing
(even rapidly increasing) sequence of functions
of the form $\sum_{i=1}^n \alpha_i\widehat{[c_i]}$, where
$\alpha_i\in \R^+$ and $[c_i]\in \Cu(A)$ for  $i=1,\dots,n$.
\end{corollary}
\begin{proof}
By \eqref{ahat}---cf.~Proposition \ref{quasifunct}---, the function 
$\widehat a$ is the supremum of the
increasing sequence 
of Riemann sums $\sum_{i=1}^{2^n} 1/2^n\widehat{[(a-i/2^n)]}$. Let $f$
be an arbitrary
function in $\LL(\F(\Cu(A)))$. Then by Theorem \ref{supsofas}
there is a rapidly increasing sequence
$(\widehat{a_i})$,
with $a_i\in (A\otimes \mathcal K)^+$, such that $f=\sup \widehat{a_i}$.
We may now
interpolate between
$\widehat{a_i}$ and $\widehat{a_{i+i}}$ an element of the form $\sum_{i=1}^n \alpha_i\widehat{[c_i]}$.
This proves the result.
\end{proof}

\begin{corollary}
For every functional  $\Lambda\colon \LL(\F(\Cu(A)))\to [0,\infty]$ (i.e,
additive, order-preserving
map, taking 0 into 0, and preserving suprema of increasing sequences)
there is $\lambda\in \F(\Cu(A))$
such that $\Lambda(f)=f(\lambda)$.
\end{corollary}  
\begin{proof}
Set $\Lambda(\widehat{[a]})=\lambda([a])$. Then $\Lambda(f)=f(\lambda)$
for every $f$ of the form $\widehat{[a]}$.
By the previous corollary this equality also holds all $f\in \LL(\F(\Cu(A)))$.
\end{proof}

Let $\phi\colon A\to B$ be homomorphism of C*-algebras. Recall that
$\F(\phi)$ is a continuous linear
map from $\F(\Cu(B))$ to $\F(\Cu(A))$. It follows that
$f\mapsto f\circ \F(\phi)$ maps $\LS(\F(\Cu(A)))$ to 
$\LS(\F(\Cu(B)))$.
Moreover, in this way $\widehat a$ is mapped to $\widehat {\phi(a)}$.
Thus, by Theorem \ref{supsofas}, 
$\LL(\F(\Cu(\phi)))(f):=f\circ \F(\Cu(\phi))$ is a map from
$\LL(\F(\Cu(A)))$ to $\LL(\F(\Cu(A)))$. 

\begin{theorem}
$\LL(\F(\Cu(\cdot)))$ is a sequentially continuous covariant functor from 
the category of C*-algebras to the category $\mcCu$. 
\end{theorem}
\begin{proof}
Let us first show that $\LL(\F(\Cu(A)))$ is an ordered semigroup in
the category $\mcCu$.
We have already seen that the supremum of an increasing sequence in
$\LL(\F(\Cu(A)))$, with respect to the 
pointwise order under consideration, exists, and is equal to the
pointwise supremum of the sequence.
By Theorem \ref{supsofas} every element is the supremum of a rapidly
increasing sequence (i.e., a sequence satisfying the axiom (2)
of the category $\mcCu$)
of functions that also belong to $\LL(\F(\Cu(A)))$. We clearly have the axiom (3)
of the category $\mcCu$ too, since
the supremum of an increasing sequence of functions in $\LL(\F(\Cu(A)))$
is the pointwise supremum of the sequence. 
Suppose that $f_1\wayb g_1$ and $f_2\wayb g_2$ in $\LL(\F(\Cu(A)))$.
Let $h_1$ be such that 
$f_1\leq h_1\leq (1-\mu)g_1$,  
and $h_1$ is continuous at each point where $g_1$ is finite. Suppose
that $h_2$ is in the same relationship
with respect to $f_2$ and $g_2$. Then
$f_1+f_2\leq h_1+h_2\leq (1-\mu)(g_1+g_2)$, and $h_1+h_2$ is continuous
at the points where
$g_1+g_2$ is finite. Hence, $f_1+f_2\wayb g_1+g_2$. This shows that
$\LL(\F(\Cu(A)))$ is in $\mcCu$.

If $\phi\colon A\to B$ is a homomorphism of C*-algebras then $\F(\Cu(\phi))$
is continuous and linear.
Keeping this in mind, it is easy to show that $\LL(\F(\Cu(\phi)))$
preserves suprema of increasing sequences and 
the relation $\wayb$. 

Let $A=\varinjlim (A_i,\phi_{i,j},i,j\in \N)$ be a sequential inductive
limit of C*-algebras. 
Recall from the proof of Theorem \ref{Dcompact} the two conditions
(1) and (2) that characterize
inductive limits in the category $\mcCu$. 
In order to show that $\LL(\F(\Cu(A)))$ is the inductive limit of
the $\LL(\F(\Cu(A_i)))$s 
it is enough to show that these conditions are satisfied with respect to 
$\LL(\F(\Cu(A)))$ and the inductive 
system $(\LL(\F(\Cu(A_i))), \LL(\F(\Cu(\phi_{i,j})))$.

Let us show that the condition (1) is satisfied. It was shown in
\cite[Theorem 2]{coward-elliott-ivanescu} that 
$\Cu(A)=\varinjlim (\Cu(A_i),\Cu(\phi_{i,j}),i,j\in \N)$ in the
category $\mcCu$. Thus, for every $[c]$ in $\Cu(A)$ there is an 
increasing sequence $([a_n])$ of elements coming from the finite stages
of the limit and with supremum $[c]$. 
It follows that the same is true for every element of 
$\LL(\F(\Cu(A)))$ of
the form $\sum_{i=1}^m \alpha_i\widehat{[c_i]}$.
Finally, since, by Corollary \ref{supsof[a]}, every $f\in \LL(\F(\Cu(A)))$
is the supremum of an
increasing---even rapidly increasing---sequence 
of functions of the form $\sum_{i=1}^m \alpha_i\widehat{[c_i]}$, we
deduce by a standard argument that $f$ is also
the supremum of an increasing sequence of elements coming from finite stages.

Let us now show that the condition (2) of the proof of Theorem
\ref{Dcompact} is satisfied. 
For $h\in \LL(\F(\Cu(A_1)))$
let us denote $\LL(\F(\Cu(\phi_{1,i})))(h)$ by $h_i$, for $i=2,\dots,\infty$. 
Let $f',f,g\in \LL(\F(\Cu(A_1)))$ be such that
$f'\wayb f$ and $f_\infty\leq g_\infty$.
Then the compact sets $\overline{\Set(f_i')}\cap \Set(g_i)^c$
have as projective limit the set $\overline{\Set(f_\infty')}\cap \Set(g_{\infty})^c$.
This last set is empty, since
$f_\infty'\wayb g_\infty$.  Therefore, for some $i$ we must have
$\overline{\Set(f_i')}\subseteq \Set(g_i)$, and so 
$f_i'\leq g_i$. 
\end{proof}

The following lemma will be used in the next section.
\begin{lemma}\label{lostlemma}
 If $[a]\wayb [b]$ then $\widehat{[a]}\wayb (1+\delta)\widehat{[b]}$
 for all $\delta>0$.
\end{lemma}

\begin{proof}
By Proposition \ref{closedbysups} (ii), it is sufficient to  show that $\overline{\Set(\widehat{[a]})}\subseteq \Set((1+\delta)\widehat{[b]})$.
Let $(\lambda_i)$ be a net in $\Set(\widehat{[a]})$ converging to $\lambda$.
By the definition
of the topology of $\F(\Cu(A))$ we have $1\leq \limsup \lambda_i[a]\leq \lambda([b])$.
So $\lambda\in \Set((1+\delta)\widehat{[b]})$ for any $\delta>0$.
\end{proof}

\subsection{The space $\LL(\T(A))$}
Here we briefly review the properties of the space $\LL(\T(A))$.

If $a\in A^+$ then 
$\bar a(\tau)=\tau(a)$ defines a lower semicontinuous function 
in $\LL(\T(A))$. All the propositions and lemmas that were proved before 
for the ordered cone $\LL(\F(\Cu(A)))$ have obvious counterparts for $\LL(\T(A))$.
The proofs of these results
are entirely analogous to the ones that we have seen above. 
We therefore  have the following theorem:
\begin{theorem}\label{supsofasT}
Let  $f$ be in $\LL(\T(A))$. Then $f$ is the supremum of an increasing
sequence $(\bar a_n)$, with $a_n\in A^+$.
\end{theorem} 
\begin{remark} Notice that the positive elements $a_n$ are now chosen in
the C*-algebra $A$ and not in $A\otimes \mathcal K$
(unlike in Theorem \ref{supsofas}). The stability of $A\otimes \mathcal K$
was used in the proof
of Theorem \ref{supsofas} to find orthogonal elements $a$ and $b$  that
were Murray-von Neumann
equivalent to two given elements. This step is not needed in proving
Theorem \ref{supsofasT},
since the additivity on pairs of orthogonal elements of quasitraces is now
replaced by the full additivity of traces.
\end{remark}

\begin{remark} 
As was done before for $\LL(\F(\Cu(A)))$, Theorem \ref{supsofasT}
may be used to show that $\LL(\T(A))$ is an ordered semigroup in the 
category $\mcCu$ and
$\LL(\T(\cdot))$ is a continuous functor from the category of C*-algebras
to the category $\mcCu$.
Theorem \ref{supsofasT} is also used, and at the same time improved, in
\cite[Theorem 2]{tracecomp}. It is shown
there that if $A$ is stable then $f\in \LL(\T(A))$ if and only if
$f=\bar a$ for some $a\in A^+$.
\end{remark}

\section{The structure of the Cuntz semigroups of certain C*-algebras}
Let us apply the results of the previous sections to study the structure of
the Cuntz semigroup, in certain well-behaved cases.

\begin{proposition}
Suppose that $A$ is an AH C*-algebra with no dimension growth. Then there 
is a constant $M$ such that 
for all $[a]$ and $[b]$ in $\Cu(A)$ 
we have $\widehat{[a]}\leq \widehat {[b]}$ if and only if
$k[a]\leq (k+M)[b]$ for all $k\in \N$.
\end{proposition} 

\begin{proof}
First suppose that $A$ is a direct sum of homogeneous algebras. 
For $a\in (A\otimes \mathcal K)^+$ let $\rank([a])$ denote the lower
semicontinuous function
on the spectrum of $A$ given by $\rank([a])(x):=\rank(a(x))$.
In this case we shall obtain the desired result as a corollary of the
following theorem of Toms (see \cite[Theorem 3.15]{toms}): 

\emph{There is a constant $K$
such that for every finite dimensional  compact Hausdorff
space $X$, if $[a],[b]\in
\Cu(C_0(X))$ satisfy  
$\rank[a]+K\dim X\leq \rank [b]$, then $[a]\leq [b]$.}
 
Let us see how.
Suppose that $\widehat{[a]}\leq \widehat{[b]}$.
This implies that $\rank[a]\leq \rank [b]$. So 
$k \cdot\rank[a]+K\dim X\leq k \cdot\rank [b]+K\dim X$
for all $k\in \N$. We may assume without loss of generality that
$b(x)\neq 0$ for all $x$. So $\rank b\geq 1$
and $k \cdot\rank[a]+K\dim X\leq (k +K\dim X)\cdot\rank [b]$. We conclude
that $k[a]\leq (k+M)[b]$
for $M=K\dim X$ and all $k\in \N$.

Let $A$ be an AH algebra with no dimension growth. Suppose that
$A=\varinjlim (A_i,\phi_{i,j},i,j\in \N)$
where the $A_i$s are homogeneous algebras with spectra  of bounded dimension.
Since $\LL(\F(\Cu(\cdot)))$ is a sequentially continuous functor,
$\LL(\F(\Cu(A)))$ is the limit of the $\LL(\F(\Cu(A_i)))$s in 
the category $\mcCu$.
First let us suppose that $[a]$ and $[b]$ come from finite stages
of the sequence $(A_i)$, say, $[a]=\phi_{1,\infty}([a_1])$ and
$[b]=\phi_{1,\infty}([b_1])$.
Let us write $a_n=\phi_{1,n}(a_1)$ and $b_n=\phi_{1,n}(b_1)$.
For every $\epsilon>0$ and $k>0$ we have 
$\widehat{[(a_1-\epsilon)_+]}\wayb (1+1/k)\widehat{[a_1]}$
by Lemma \ref{lostlemma}. By the condition (2) for  inductive limits in the 
category $\mcCu$ (see proof of Theorem \ref{Dcompact}), there is $n$ 
such that
$\widehat {[(a_n-\epsilon)_+]}\leq (1+1/k)\widehat {[b_n]}$. Therefore,
$k\widehat {[(a_n-\epsilon)_+]}\leq (k+1) \widehat {[b_n]}$. We have
already established that as this is at a finite stage
of the sequence $(A_i)$ this implies that $k[(a_n-\epsilon)_+]\leq (k+M) [b_n]$,
where $M$ depends only on the bound on the dimensions of the spectra of
this $A_i s$---and, as we may take the best choice for different
inductive limit decompositions, therefore depends only on $A$.
Recalling that $a_n$ and $b_n$ map into $[a]$ and $[b]$ in $A$,
we have $k[(a-\epsilon)_+]\leq (k+M) [b]$.

Let us consider the case that $[a]$ and $[b]$ do not necessarily come 
from finite stages. Let $[f_i]$ and $[g_i]$, $i=1,2,\dots$, be increasing
sequences of
elements coming from finite stages and with suprema $[a]$ and $[b]$ respectively.
Since $\widehat {[f_i]}\wayb \widehat {[a]}\leq \widehat {[b]}$,
there is $[g_j]$ such that $\widehat {[f_i]}\leq \widehat {[g_j]}$. Hence as
before, $k[f_i]\leq (k+C)[g_j]\leq (k+M)[b]$.
Passing to the supremum with respect to $i$ we get that $k[a]\leq (k+M)[b]$.
\end{proof}

We now turn to the setting of arbitrary C*-algebras with almost unperforated
Cuntz semigroup. Recall from \cite{rordam} that an ordered semigroup 
is said to be almost unperforated if the inequality $(k+1)x\leq ky$ for some
$k\in \N$ implies that $x\leq y$.

The following proposition is an improvement of 
\cite[Proposition 3.2]{rordam} for ordered semigroups 
in the category $\mcCu$.

\begin{proposition}\label{improverordam}
Let $S$ be an ordered semigroup in the category $\mcCu$. 
Then $S$ is almost unperforated 
if and only if the following condition is fulfilled: for all $x$ and $y$ in $S$ with $x\leq \infty \cdot y$ and 
$\lambda(x)<\lambda(y)$ for any functional on $S$ such that $\lambda(y)=1$, one has $x\leq y$.
\end{proposition}

\begin{proof}
Suppose that $S$ satisfies the condition of comparison of elements by
functionals described in the
statement of the proposition (this condition is often referred to as
``strict comparison''). If $(k+1)x\leq ky$
then $\lambda(x)\leq k/(k+1)<\lambda(y)$ for any $\lambda$ such that $\lambda(y)=1$. Since $x\leq ky\leq \infty\cdot y$,
we conclude that $x\leq y$, as desired.

Suppose that $S$ is almost unperforated. Let $x,y\in S$ be such that $x\leq \infty \cdot y$ and $\lambda(x)<\lambda(y)$
for all $\lambda$ such that $0<\lambda(y)<\infty$. Let $z\wayb x$. Then $z\leq ky$ for some $k$.    
We shall prove that for every additive, order preserving, function $D$ on $S$---not necessarily preserving suprema of 
increasing sequences---such that $D(y)=1$, we have $D(z)<D(y)$. By \cite[Proposition 3.2]{rordam}, this will imply 
that $z\leq y$, from which the desired result will follow on taking the supremum over all $z$ that are far below $x$.

Let $D\colon S\to [0,\infty]$ be additive, order preserving, and such that $D(y)=1$. Set 
$\widetilde Dw:=\sup\{\, Dw'\mid w'\wayb w\,\}$, for each $w\in S$. By Lemma \ref{whatsaD}, $\widetilde D$ is a 
functional on $S$, and it is clear that $\widetilde D(y)\leq 1<\infty$.

Case 1. Suppose that $\widetilde Dy\neq 0$. Then 
\[
D(z)\leq \widetilde D(x)<\widetilde D(y)\leq D(y).
\]  

Case 2. Suppose that $\widetilde Dy= 0$. Then $\widetilde Dx=0$ 
(because $x\leq \infty\cdot y$), and so
\[
D(z)\leq \widetilde D(x) =0 < D(y).
\qedhere
\]
\end{proof}

\begin{corollary}\label{semigroupHB}
Let $S$ be an almost unperforated ordered semigroup in the 
category $\mcCu$. 
Let $x,y\in S$. Then $\lambda(x)\leq \lambda(y)$ for every functional $\lambda$ on $S$ if and only if $kx\leq (k+1)y$ for all $k\in \N$.
\end{corollary}

\begin{proof}
The implication that $\lambda(x)\leq \lambda(y)$ for every functional $\lambda$ if $kx\leq (k+1)y$ for all $k\in \N$
is obvious. Let us prove the converse.

By considering functionals with the only possible values $0$ or $\infty$ we conclude 
that $x\leq \infty y$. We may now apply the previous proposition with $kx$ and $(k+1)y$ in place
of $x$ and $y$, respectively.
\end{proof}

By a theorem of R\o rdam (see \cite[Theorem 4.5]{rordam}), the Cuntz semigroup of a C*-algebra that absorbs the Jiang-Su algebra is almost
unperforated. It follows by Corollary \ref{semigroupHB} that if $A$ absorbs the Jiang-Su algebra then 
$\widehat{[a]}\leq \widehat{[b]}$
if and only if $k[a]\leq (k+1)[b]$ for all $k\in \N$. In  the sequel we shall denote the Jiang-Su algebra
by the letter $Z$.

Let us now try to identify $\LL(\F(\Cu(A)))$ with a subsemigroup of $\Cu(A)$. Notice that if $f\in \LL(\F(\Cu(A)))$ 
is compact (i.e., $f\wayb f$ in $\LL(\F(\Cu(A)))$)
then $(1-\epsilon)f=f$ for some $0\not=\epsilon \not= 1$  
and so $f$ takes 0 and $\infty$ as its only possible values.
This observation suggests the following definition.

Let us say that $[a]\in \Cu(A)$ is purely non-compact if it has the property 
that if the image of $[a]$ in the quotient by some ideal $I$---let us denote this by 
$[a_I]$---is compact (i.e., $[a_I]\wayb [a_I]$), then $[a_I]$ is a multiple of infinity (i.e., $2[a_I]=[a_I]$).

\begin{proposition}\label{ppositive}
(i) The purely non-compact elements of $\Cu(A)$ form a subobject of $\Cu(A)$, i.e., 
a subsemigroup closed under the passage to suprema of increasing sequences.

(ii) If $[a]$ is purely non-compact then for all $[b]$, with $[b]\wayb [a]$, and for all $M\in \N$, we have
$(k+M)[b]\leq k[a]$ for all sufficiently large $k$.

(iii) If $\Cu(A)$ is almost unperforated the converse of (ii) is true.

(iv) If $A$ is simple then every element of $\Cu(A)$ is purely non-compact except for the element
$[p]$ where $p$ is any non-zero finite projection.

(v) If $[a]=\sum_{i=1}^\infty [c_i]$ and $\Ideal(c_i)=\Ideal(a)$, then $[a]$ is purely non-compact. 
\end{proposition}

\begin{proof}
Before proving the proposition we need
some preliminary formulas.

Suppose that $[a]$ is purely non-compact.
Let $\epsilon>0$ and choose a positive function $c_{\epsilon}(t)$ different from 0 precisely 
on the interval $(0,\epsilon)$. Then $(a-\epsilon)_+$ and $c_\epsilon(a)$ are orthogonal
and $(a-\epsilon)_+ + c_\epsilon(a)\leq C_1a\leq C_2((a-\epsilon/2)_+ + c_\epsilon(a))$ for some
positive scalars $C_1$ and $C_2$. Hence,
\begin{equation}\label{acepsilon}
[(a-\epsilon)_+]+[c_\epsilon(a)]\leq [a]\leq [(a-\epsilon/2)_+]+[c_\epsilon(a)].
\end{equation}
These inequalities imply that $[a]$ is compact after passing to the quotient by the ideal
$\Ideal(c_\epsilon(a))$. Let us call this ideal $I$ and
let us denote with the subscript $I$
the images of elements of $A$ in $A/I$. We have that $2[a_I]=[((a-\epsilon)_+)_I]$. By \cite[Theorem 1]{crs},
this means that $2[a]\leq [(a-\epsilon)_+]+[g]$ for some $[g]$ such 
that $\Ideal(g)=I$. 
Since $[(a-\epsilon)_+]\wayb [a]$ and $[g]$ is the supremum of an
increasing sequence $[g_i]$ with 
$[g_i]\wayb [g]$, so that also $2[(a-\epsilon)_+]\wayb 2[a]\leq 
[(a-\epsilon)_+]+[g]=\sup_i ([(a-\epsilon)_+]+[g_i])$,
there is $[g']\wayb [g]$
such that $2[(a-\epsilon)_+]\leq [(a-\epsilon)_+]+[g']$. Since $[g']\wayb \infty [c_\epsilon(a)]$, there is $k\in \N$ such that
$[g']\leq k[c_\epsilon(a)]$. Thus,
\begin{equation}\label{twoisone}
2[(a-\epsilon)_+]\leq [(a-\epsilon)_+]+ k[c_\epsilon(a)]
\end{equation}
for sufficiently large $k$.

(i) Suppose that $[a_i]$ is an increasing sequence of purely non-compact 
elements with supremum $[a]$, and that
$[a]$ is compact in $\Cu(A/I)$ for some ideal $I$. Without loss of
generality, in order to prove that $[a_I]$ is
a multiple of infinity, we may assume that $I=0$. Thus,
$[a]$ is compact, and so $[a]=[a_{i_0}]$ for some $i_0$.
Since $[a_{i_0}]$ is purely non-compact, $[a_{i_0}]$ is a multiple of infinity,
and therefore so also is $[a]$, as desired.

Let $[a]$ and $[b]$ be purely non-compact and suppose that $[a]+[b]$
is compact in $\Cu(A/I)$ for some ideal
$I$. Again we may assume that $I=0$. Thus, $[a]+[b]$ is compact, and so
\[
[(a-\epsilon)_+]+[(b-\epsilon)_+]+[c_\epsilon(a)]+[c_\epsilon(b)]\leq [a]+[b]=[(a-\epsilon)_+]+[(b-\epsilon)_+]
\]
for some $\epsilon$. Let $k$ be such that \eqref{twoisone} holds
both for $[a]$ and for $[b]$. Then we have
\[
[a]+[b]=[(a-\epsilon)_+]+[(b-\epsilon)_+]+k([c_\epsilon(a)]+[c_\epsilon(b)])\geq 2([a]+[b]).
\]
Therefore, $[a]+[b]$ is a multiple of infinity, as desired.

(ii) It is enough to assume that $[b]=[(a-\epsilon)_+]$. From \eqref{twoisone} we get by
induction that for every $M\in \N$ we have $M[(a-\epsilon)_+]\leq [(a-\epsilon)_+]+ k[c_\epsilon(a)]$
for sufficiently large $k$. On the other hand we deduce from \eqref{acepsilon} that
$k[(a-\epsilon)_+]+k[c_\epsilon(a)]\leq k[a]$. Combining these two equations we get the desired
result.

(iii) Suppose that $\Cu(A)$ is almost unperforated and let $[a]$ be such that for every $\epsilon>0$ we have $k[(a-\epsilon)_+]\leq (k+1)[a]$ for $k$ large enough.
Suppose that $[a]$ is compact. Then for some $\epsilon>0$ we have $[(a-\epsilon)_+]=[a]$. So $(k+1)[a]=k[a]$,
and this implies that $2(k+1)[a]=k[a]$. Since $\Cu(A)$ is almost unperforated
it follows that $2[a]=[a]$. The same implication holds in any 
quotient of $A$.

(iv) Suppose that $A$ is simple. Then for any 
$0\not= a \in (A\otimes \mathcal K)^+$ the element $\infty\cdot [a]$ of 
$\Cu (A)$ is the same, and is the largest
element of $\Cu (A)$. (In general, the Cuntz 
semigroup may not have a largest element.)

Fix a $\in(A\otimes \mathcal K)^+$. The element $[a]$ of $\Cu(A)$
is purely non-compact unless $[a]$ is compact and not equal either to 0 
or to the largest element of $\Cu(A)$.
Suppose that the latter is true, i.e., that $[a]$ is compact
and different from both zero and the largest element of $\Cu(A)$.
For some $\epsilon>0$ we have
$[(a-\epsilon)_+]=[a]$. By \eqref{acepsilon} we have
$[a]+[c_\epsilon(a)]=[a]$. If $c_\epsilon(a)$
were non-zero then, again as $A$ is simple,  $\infty\cdot[c_\epsilon(a)]$
would be  the largest element of $\Cu(A)$, and hence 
$[a]$ would be the largest element of $\Cu(A)$. 
Therefore, $c_\epsilon(a)=0$, and  so
$[a]=[p]$ for some non-zero projection $p$.
The projection $p$ must be finite, as otherwise $[p]$ would be
the largest element of $\Cu(A)$.

(v) Suppose that $[a]$ is compact.
Then $[a]=[s]+[r]=[s]$ for some $[r],[s]\in \Cu(A)$ such that
$\Ideal(r)=\Ideal(s)=\Ideal(a)$ 
($[s]$ may be chosen as a partial sum in the given representation of $[a]$ 
as an infinite series and $[r]$ as the remainder term). Since
$\Ideal(r)=\Ideal(a)$ we have $[a]\leq \infty [r]$. Since $[a]$ is compact
we therefore have $[a]\leq k[r]$ for some $k\in \N$. We now have
$[a]=[s]+k[r]\geq 2[a]$.
The same implication holds in any quotient of $A$.
\end{proof}

\begin{lemma}\label{thereisc}
If $A$ absorbs the Jiang-Su algebra then for every $a\in (A\otimes \mathcal K)^+$ there is $[c]\in \Cu(A)$ such that
$\widehat a=\widehat{[c]}$. One can always choose $[c]$ to be purely non-compact.
\end{lemma}

\begin{proof}
First let us prove that if the existence of $[c]$ is guaranteed (in general), then we can always choose it
so that it is purely non-compact. Let $[c_i]$, $i=1,2,\dots$, be such that $\widehat a/2^i=\widehat{[c_i]}$ for all $i$. Then $[c]=\sum_{i=1}^\infty [c_i]$ is purely non-compact by Proposition \ref{ppositive} (v), 
and $\widehat{[c]}=\sum_{i=1}^\infty \widehat{[c_i]}=\widehat{a}$.

Let us now prove that $[c]$ exists with $\widehat{[c]}=\widehat a$.
Every positive element of $A\otimes \mathcal K\otimes Z$  is approximately unitarily equivalent to one of the form $b\otimes 1$
with $b\in (A\otimes \mathcal K)^+$ (see the proof of \cite[Theorem 5.5]{brown-perera-toms}). Therefore, we may assume that the given positive element 
$a$ has the form  $b\otimes 1\in (A\otimes \mathcal K)\otimes Z$.
Recall that 
\[
\widehat {(b\otimes 1)}(\lambda)=\int_{0}^\infty \lambda([(b-t)_+\otimes 1])\,dt.
\] 
Since $t\mapsto \lambda([(b-t)_+])$ is a decreasing function
of $t$, the Riemann sums
\[
\sum_{i=1}^{2^n} \frac{1}{2^n}\lambda([(b-\frac {i}{2^n})_+\otimes 1])
\]
converge to the integral. Choose a positive element $e_{1/2^n}$ of the Jiang-Su
algebra with rank $1/2^n$. (That such an element exists follows, for instance, by the computation of the
Cuntz semigroup of the Jiang-Su algebra obtained in \cite{brown-perera-toms}.) 
Then the Riemann sum above is equal to
\[
\lambda\left(\sum_{i=1}^{2^n} [(b-\frac {i}{2^n})_+\otimes e_{1/2^n}]\right).
\]
Let us show that the Cuntz semigroup elements 
\begin{align}\label{riemann}
\sum_{i=1}^{2^n} [(b-\frac {i}{2^n})_+\otimes e_{1/2^n}]
\end{align} 
form an increasing sequence. Comparing two consecutive terms of this
sequence we see that it is enough to show that
$[(b-i/2^n)_+\otimes e_{1/2^n}]=2[(b-i/2^n)_+\otimes e_{1/2^{n+1}}]$.
This is true, since $[e_{1/2^n}]=2[e_{1/2^{n+1}}]$, as follows from
the computation of $\Cu(Z)$ in \cite{brown-perera-toms}.  We have
$\widehat {b\otimes 1}=\widehat{[c]}$, with 
$[c]$ the supremum of the sequence \eqref{riemann}.  
\end{proof}

\begin{theorem}\label{mainnoncompact} 
Let $A$ be a C*-algebra.
Suppose that $\Cu(A)$ is almost unperforated. If $[a]$ and $[b]$ are in $\Cu(A)$ and $[a]$
is purely non-compact, then $[a]\leq [b]$ if and only if $\widehat{[a]}\leq \widehat {[b]}$. 

Suppose further that $A$ absorbs the Jiang-Su algebra. Then 
the map $[a]\mapsto \widehat{[a]}$ is an isomorphism of the ordered semigroups of 
purely non-compact elements of $\Cu(A)$ and of $\LL(\F(\Cu(A)))$.
\end{theorem}

\begin{proof}
Suppose that $[a]$ is purely non-compact, $[b]\in \Cu(A)$, and
$\widehat{[a]}\leq \widehat {[b]}$. By Corollary \ref{semigroupHB}, for every $k\in \N$ we have $k[a]\leq (k+1)[b]$. 
On the other hand, by Proposition \ref{ppositive} (ii), for every $\epsilon>0$  there exists $k$ such that $(k+2)[(a-\epsilon)_+]\leq k[a]\leq (k+1)[b]$. 
Since $\Cu(A)$ is almost unperforated, it follows  that $[(a-\epsilon)_+]\leq [b]$ for all $\epsilon>0$. Hence, $[a]\leq [b]$.

Let us now prove that the map $[a]\mapsto \widehat{[a]}$ is a surjection 
from the purely non-compact elements to $\LL(\F(\Cu(A)))$.
Let $f\in \LL(\F(\Cu(A)))$. By Theorem \ref{supsofas} there exists an increasing sequence $(\widehat {a_i})$
with supremum $f$, where $a_i\in (A\otimes \mathcal K)^+$. By Lemma \ref{thereisc}, for each $i=1,2,\dots$ there exists
a purely non-compact element $[c_i]\in \Cu(A)$ such that $\widehat {a_i}=\widehat{[c_i]}$. We have $\widehat{[c_i]}\leq \widehat{[c_{i+1}]}$,
and so as shown above $[c_i]\leq [c_{i+1}]$; that is, $[c_i]$, $i=1,\dots$, is an increasing sequence. 
Set $\sup [c_i]=[c]$. Since the map $[a]\mapsto \widehat{[a]}$ preserves suprema of increasing sequences, we have $f=\widehat{[c]}$. In order to ensure that $[c]$ is purely
non-compact let us choose $[c_i']$ such that $1/2^i f=\widehat{[c_i']}$ for all $i=1,2,\dots$. Then $[c']=\sum_{i=1}^\infty [c_i']$ is purely non-compact
by Proposition \ref{ppositive} (i), and  $f=\widehat{[c']}$.
\end{proof}

\begin{corollary}\label{nosimpleq}
Let $A$ be a C*-algebra that absorbs the Jiang-Su algebra.
If $A$ has no non-zero simple subquotients then $\Cu(A)\cong \LL(\F(\Cu(A)))$.
(Here by $\cong$ is meant a natural isomorphism of functors into the
category $\mcCu$.)
\end{corollary}

\begin{proof}
If $A$ has no non-zero simple subquotients then  no non-zero element of 
$\Cu(A)$, or of $\Cu(A/I)$ for an ideal $I$,
is compact (and in particular, every element of $\Cu(A)$ is purely
non-compact, as desired). For if $[a]$ is compact  and non-zero then
$[a]=[(a-\epsilon)_+]$ for some 
$\epsilon>0$, whence $\Ideal(a)=\Ideal((a-\epsilon)_+)$.
It follows, for instance from Lemma \ref{K-R}, that $\Ideal(a)$ is
a compact ideal of $A$
(a compact element of $\lat(A)$).
(Here we do not mean only countable compactness; however,
that would be sufficient for our purposes since $\Ideal(a)$
is singly generated.)
Hence, $\Ideal(a)/J$ is simple
for any maximal proper ideal $J$
of $\Ideal(a)$ ($J$ exists, by Zorn's Lemma, since $\Ideal(a)$ is
compact and non-zero). 
\end{proof}

As another corollary of Theorem \ref{mainnoncompact}, let us give a 
computation of the Cuntz
semigroup of a simple C*-algebra absorbing the Jiang-Su algebra
(this computation was previously
obtained in \cite{brown-perera-toms} with the additional assumptions
that the algebra 
was unital, exact, and of stable rank one). Let $A$ be a simple
C*-algebra absorbing
the Jiang-Su algebra. Let $\mathrm{V}(A)$ denote the semigroup of
Murray-von Neumann equivalence  classes of projections
of $A\otimes \mathcal K$.  Let us define on the abstract disjoint union 
$(\mathrm{V}(A)\backslash \{0\})\sqcup \LL(\F(\Cu(A)))$ 
an order and an addition operation, making it what might be 
called the lexicographic ordered semigroup disjoint union. 
Inside the two subsets 
$(\mathrm{V}(A)\backslash \{0\})$ and $\LL(\F(\Cu(A)))$ 
let us retain the order and addition with which these sets are endowed.
Let $[p]\in \mathrm{V}(A)\backslash\{0\}$ 
and $f\in \LL(\F(\Cu(A)))$. Let us define $f+[p]\in \LL(\F(\Cu(A)))$
as the
function $f+\widehat{[p]}$ if $f\not= 0$, and as the class $[p]$ if
$f=0$. Let us say that   
$f\leq [p]$ if $f\leq \widehat{[p]}$, and that $[p]\leq f$ if
$\widehat{[p]}+ g=f$ for some $0\not= g\in L (F(\Cu(A)))$. It
is not difficult to verify that $(\mathrm{V} (A)\backslash \{0\})
\sqcup \LL(\F (\Cu(A)))$ is an object in $\mcCu$.

\begin{corollary}
Let $A$ be a simple C*-algebra absorbing
the Jiang-Su algebra. Then  either $A$ is purely infinite and
$\Cu(A)\cong \{0,\infty\}$,
or every projection in $A\otimes \mathcal K$ is finite and
$\Cu(A)\cong  (\mathrm{V}(A)\backslash\{0\})\sqcup \LL(\F(\Cu(A)))$.
(Here by $\cong$ is meant a natural isomorphism of 
functors into the category $\mcCu$.)
\end{corollary}

\begin{proof}
If $\F(\Cu(A))\cong \{0,\infty\}$ then by Proposition \ref{improverordam},
$\Cu(A)\cong \{0,\infty\}$, in other words, $A$ is purely infinite 
(this is also obtained in \cite[Corollary 5.1]{rordam}).  If, on the
other hand, $\F(\Cu(A))$ contains a non-trivial functional,
i.e., a functional with a non-zero finite value, then every 
projection in $A\otimes \mathcal K$ must be finite.
Indeed, if $\lambda$ is a non-trivial functional, and $p$ is
a projection in $A\otimes \mathcal K$, 
then, by \cite{coward-elliott-ivanescu}, the class $[p]$ is compact in $\Cu(A)$,
and so by simplicity is majorized by a finite multiple of any non-zero element, and
hence is finite on $\lambda$. Again by simplicity, 
Ker $\lambda =0$, and it follows that if $[p]+[a]=[p]$ in
$\Cu (A)$ then $[a]=0$, and in particular the projection $p$
is finite.

Suppose that every projection of $A\otimes \mathcal K$ is finite. 
By Proposition \ref{ppositive} (iv), the complement of the set of 
purely non-compact elements of $\Cu(A)$, in the case of a 
simple C*-algebra, 
is the set of elements $[p]$ such that $p$ is a non-zero finite projection.
It is easy to show (and well known) 
that among finite projections Cuntz equivalence amounts to
Murray-von Neumann equivalence. Therefore, by Theorem \ref{mainnoncompact},  
the map from $\Cu(A)$ to
$(\mathrm{V}(A)\backslash\{0\})\sqcup \LL(\F(\Cu(A)))$ given by
$[p]\mapsto [p]$ if $p$ is a projection, and 
$[a]\mapsto \widehat{[a]}$ if $[a]$ is purely non-compact,
is a bijection. 

To prove that that this (natural) map is an isomorphism
of ordered semigroups, let $0\not= p \in A\otimes \mathcal K$ 
be a projection and let $0\not=a \in (A\times \mathcal K)^+$
be such that $[a]$ is purely non-compact in $\Cu (A)$.
The sum $[p]+[a]$ in $\Cu(A)$ is
then also purely non-compact. (By \ref{ppositive} (iv) it is enough to
show that $[p]+[a]$ is not the class of a non-zero finite
projection, if $[a]$ itself is not the class of a non-zero
finite projection. We may assume that 
$pa=0$. If $[p]+[a]=[p+a]$ is the class of a non-zero finite
projection, say $q$, then, by the formulation of
Cuntz equivalence given in \cite{coward-elliott-ivanescu}, and since, also by \cite{coward-elliott-ivanescu},
$[q]$ is compact, and $\sup_{\epsilon>0}[(a-\epsilon)_+]=[a]$,
whence $\sup_{\epsilon>0}[p+(a-\epsilon)_+]=[p+a]$,
the Hilbert $A$-module $qA$ is (by compactness) isomorphic
to a sub Hilbert module $X$ of $((p+(a-\epsilon)_+)A)^-$ for some
$\epsilon>0$, and (for any $\epsilon >0$) $((p+(a-\epsilon)_+)A)^-$ 
is isomorphic to a sub Hilbert module $Y$ of $qA$. Hence,
if $X'$ denotes the isomorphic copy of $X$
contained in $((p+(a-\epsilon)_+)A)^-$, and $X''$ and $Y'$
the resulting isomorphic images of $X'$
and $Y$ in $qA$, so that $X''\subseteq Y'\subseteq qA$ and $X''$ is
isomorphic to $qA$, by finiteness of $q$ it follows
that $X''=qA$---so that in particular $Y=qA$, and
so $((p+(a-\epsilon)_+)A)^-=q'A$ for some projection 
$q'=q'_\epsilon$. Then
necessarily $p+(a-\epsilon)_+=q'$, i.e., $(a-\epsilon)_+$ is a
projection (namely, $q'-p$). Since $\epsilon>0$ may be
arbitrarily small, it follows that $a$ itself is a projection,
equal to $(a- \epsilon)_+$ for	some
$\epsilon>0$, and therefore $\le q'_\epsilon$ for that $\epsilon$;
and therefore finite.
By hypothesis, $[a]\not= 0$; we have therefore proved the
contrapositive.)

It follows that the image of $[p]+[a]$ in
the disjoint union is $([p]+[a])^{\wedge}$. On
the other hand, the image of $[p]$ is $[p]$ and the 
image of $[a]$ is $\widehat{[a]}$, and the
sum of $[p]$ and $\widehat{[a]}$ in the (lexicographic) disjoint 
union is by definition $\widehat{[p]}+\widehat{[a]}$, which is equal to
$([p]+[a])^{\wedge}$. The one case remaining in which to
check additivity, that the first element is $[p]$ and
the second is the purely non-compact element $0\in \Cu A$,
is trivial: the
sum of these elements is $[p]$, and the sum of the
images, $\widehat{[p]}$ and $0$, is $\widehat{[p]}$, the image of the sum.

It remains to verify that the relation $[a]\le [p]$ or
$[p]\le [a]$ in $\Cu (A)$ implies the relation
$\widehat{[a]}\le [p]$ or $[p]\le \widehat{[a]}$ for the
images in the disjoint union, and conversely (with
$p$ and $a$ still as above).
If $[a]\le [p]$ in $\Cu (A)$ then $\widehat{[a]}\le \widehat{[p]}$ in
$L(F(\Cu (A)))$, which is just the definition of 
$\widehat{[a]}\le [p]$ in the disjoint
union. Conversely, if $\widehat{[a]}\le [p]$ in the
disjoint union, i.e., if $\widehat{[a]}\le \widehat{[p]}$ in 
$L(F(\Cu (A))$, then, as $a$ is purely non-compact,
by Theorem \ref{mainnoncompact}, $[a]\le [p]$. 

If $[p]\le [a]$ in $\Cu (A)$, then
the compact Hilbert module $pA$ is isomorphic (by \cite{coward-elliott-ivanescu})
to a sub Hilbert module of $(aA)^-$, necessarily
complemented: 
$(aA)^- \cong pA\oplus (bA)^-$ (with, e.g., $b=(1-p)a)$.
Since $a$ is non-zero and purely non-compact, $b\not= 0$.
Hence, $\widehat{[a]}=\widehat{[p]}+g$ with $0\not= 
g=\widehat{\, [b]\, }\in L
(F(\Cu (A)))$, i.e.,
$[p]\le \widehat{[a]}$ in the disjoint union. Conversely, if
$\widehat{[a]}=\widehat{[p]}+g$ with $0\not= g\in L (F(\Cu (A)))$, then
by Theorem \ref{mainnoncompact} there exists a purely non-compact element
$0\not=[b]\in \Cu (A)$ such that $\widehat{\, [b]\, }=g$. Then,
as shown above, also $[p]+[b]$ is
purely non-compact, and since $\widehat{[a]}=\widehat{[p]}
+{\widehat{\, [b]\, }}=
([p]+[b])^{\wedge}$, and also $[a]$ is non-compact, by
Theorem \ref{mainnoncompact}, $[a]=[p]+[q]$, and in particular
$[p]\le [a]$ in $\Cu (A)$.
\end{proof}

\begin{remark}
Certain C*-algebras with no non-zero simple subquotients
were considered by the second author 
in \cite{robert}---these were closed two-sided ideals of AI
algebras---and were classified by means
of tracial data. (The techniques of \cite{robert} apply in fact to
arbitrary ideals of AI algebras
with no non-zero simple subquotients---equivalently, with the
K$_1$-group of every ideal equal to zero.) For 
these algebras, the results of \cite{robert} may therefore 
be viewed as a determination 
of the Cuntz semigroup in terms of tracial data, although
the way this ordered semigroup is determined is only implicit.
Note that, in \cite{ciuperca-elliott}, the Cuntz semigroup, together with the 
special element consisting of  the class of the strictly positive
elements, was shown to be a complete invariant for arbitrary 
AI algebras, or ideals of AI algebras.
The results of \cite{robert} could be deduced from this 
together with Corollary \ref{nosimpleq}. The problem of 
describing the Cuntz semigroup
in terms of K-theoretical and tracial data
in this more general setting would seem to be very interesting. 
\end{remark}

\end{document}